\documentclass[a4paper]{article}
\usepackage{amscd,amsmath,amsthm,amssymb}
\begin{document}

\newtheorem{thm}{Theorem}[section]
\newtheorem{lem}[thm]{Lemma}
\newtheorem{prop}[thm]{Proposition}
\newtheorem{cor}[thm]{Corollary}
\newtheorem{defn}[thm]{Definition}
\newtheorem{rem}[thm]{Remark}
\newtheorem{exam}[thm]{Example}
\newtheorem{rems}[thm]{Remarks}
\newtheorem{exs}[thm]{Examples}

\newcommand{\grade}[1]{\begin{bf}#1\end{bf}}
\newcommand{\cat}[1]{\begin{bf}#1\end{bf}}
\newcommand{\mapor}[1]{\smash{\mathop{\longrightarrow}\limits^{#1}}}
\newcommand{\mapver}[1]{\Big\downarrow
\rlap{$\vcenter{\hbox{$\scriptstyle#1$}}$}}
\newcommand{\liminv}{\smash{\mathop{\lim}\limits_{\leftarrow}}}

\newcommand{\AN}{\cat{C}}
\renewcommand{\bar}{\overline}
\newcommand{\solose}{\Rightarrow}
\newcommand{\upar}[1]{\!\uparrow\!\!{#1}}
\newcommand{\downar}[1]{\!\downarrow\!\!{#1}}

\newcommand{\specif}[2]{\left\{#1\,\left|\, #2\right. \,\right\}}

 \newenvironment{pf}{\paragraph{Proof}}{\par\smallskip}
 \newenvironment{pfof}[1]{\paragraph{Proof of #1}}{\par\smallskip}
 \newenvironment{step}[1]{\subparagraph{{\sc Step #1:}}}{\par\smallskip}
\newenvironment{ex}{\begin{exam}\em}{\end{exam}}

\newcommand{\ds}{\displaystyle}
\newcommand{\oo}{\infty}
\newcommand{\de}{\partial}
\newcommand{\per}{\!\cdot\!}
\newcommand{\sA}{{\mathcal A}}
\newcommand{\sC}{{\mathcal C}}
\newcommand{\Oh}{{\mathcal O}}
\newcommand{\sM}{{\mathcal M}}
\newcommand{\sL}{{\mathcal L}}
\newcommand{\sB}{\mathcal B}

\newcommand{\Aut}{\operatorname{Aut}}
\newcommand{\Def}{\operatorname{Def}}
\newcommand{\Hom}{\operatorname{Hom}}
\newcommand{\HOM}{\operatorname{\mathcal H}om}
\newcommand{\Spec}{\operatorname{Spec}}
\newcommand{\rank}{\operatorname{rank}}
\newcommand{\Der}{\operatorname{Der}}
\newcommand{\Coder}{{\operatorname{Coder}}}

\newcommand{\somdir}[2]{\hbox{$\mathrel
{\smash{\mathop{\mathop \oplus\limits _{#1}}
\limits^{#2}}}$}}
\newcommand{\tensor}[2]{\hbox{$\mathrel
{\smash{\mathop{\mathop \otimes\limits _{#1}}
^{#2}}}$}}
\newcommand{\symm}[2]{\hbox{$\mathrel
{\smash{\mathop{\mathop \odot\limits _{#1}}
^{#2}}}$}}
\newcommand{\external}[2]{\hbox{$\mathrel
{\smash{\mathop{\mathop \bigwedge\limits_{#1}}
^{\!#2}}}$}}

\newcommand{\dual}{^{\vee}}
\newcommand{\Q}{\mathbb Q}
\newcommand{\C}{\mathbb C}
\newcommand{\A}{\mathbb A}
\newcommand{\N}{\mathbb N}
\newcommand{\Z}{\mathbb Z}
\newcommand{\K}{{\mathbb K}\,}

\title{Extended deformation functors}
\author{Marco Manetti\thanks{partially 
supported by Italian MURST program 
'Spazi di moduli e teoria delle rappresentazioni'. Member of GNSAGA 
of CNR.}\\
Universit\`a di Roma ``La Sapienza'', Italy}

\maketitle

\begin{abstract}
We introduce a precise notion, in terms of some Schlessinger's type 
conditions, of extended deformation functors which 
is compatible with most of recent ideas in the Derived Deformation 
Theory (DDT) program and with geometric examples. With this notion we 
develop the (extended) analogue of Schlessinger and obstruction 
theories. The inverse mapping theorem holds for 
natural transformations of extended deformation functors and all such 
functors with finite dimensional tangent space are prorepresentable in 
the homotopy category.  Finally we prove that the primary obstruction 
map induces a structure of graded Lie algebra on the tangent space.
\end{abstract}

			\section*{Introduction}

This paper is devoted to the study of 
some foundations in the derived deformation theory (DDT) program, 
see \cite{ciokap}, \cite{BK}. 

One of the main aspects of DDT program is the {\em Hidden Smoothness 
Philosophy}: this is based on the idea that, in characteristic 0, every 
``reasonable'' deformation problem give rise  naturally to an extended 
``quasismooth''\footnote{Some authors prefer to say ``smooth in an 
appropriate sense''}  differential $\Z$-graded moduli 
space $\sM^{\Z}$ such that the ordinary 
moduli space $\sM$ is just the truncation in degree 0 of $\sM^{\Z}$.
Roughly speaking a differential $\Z$-graded space is a 
ringed space whose structure 
sheaf takes value in differential $\Z$-graded connected algebras; see 
\cite{ciokap}, \cite{grassi2} for details.
In its first written reference \cite{K3}, the DDT program was proposed as 
a possible way to define virtual tangent bundle and virtual 
fundamental class on moduli spaces, while more recent applications of 
extended moduli spaces \cite{BK}, \cite{Baratesi}, concern the 
understanding of mirror symmetry for Calabi-Yau manifolds.

The hidden smoothness philosophy applied to infinitesimal deformation 
theory  should implies that, over a field of characteristic 0, every 
reasonable deformation functor is the truncation of a suitable 
quasismooth extended deformation functor. Formally a  
functor $F$ from the category $\cat{Art}_{\K}$ of local Artinian 
$\K$-algebras with residue field $\K$ to the category of sets is 
called a deformation functor if it satisfies the (slightly modified, 
see the discussion in \cite{FM1}, \cite{FaMa2} and references therein) 
Schlessinger's conditions 
\begin{description}
\item[1)] $F(\K)=\{0\}$.
\item[2)] $F(A\times_{\K}B)=F(A)\times F(B)$.
\item[3)] The map $F(A\times_{C}B)\to F(A)\times_{F(C)}F(B)$ is 
surjective for every pair of morphism $\alpha\colon A\to C$, 
$\beta\colon B\to C$ with $\alpha$ surjective.
\end{description}
and the term ``reasonable'' means essentially that $F$ arises from a 
geometric deformation problem.

By an extended deformation functor we understand a set valued 
functor $F$
defined on a category $\AN$ containing $\cat{Art}_{\K}$ as a full 
subcategory and such that $F$ satisfies some extended 
Schlessinger's conditions.

As far as we known there are in literature two main notions of 
extended deformation functors. 

The first (\cite{BK}, \cite{Baratesi}) 
considers $\AN$ as the category of $\Z$-graded local Artinian 
$\K$-algebras with residue field $\K$ and the extended Schlessinger's 
conditions are nothing else than the trivial extension of 1), 2) and 3). 
This approach does not present any additional difficulty with respect to 
the classical case but works well only when the 
associated extended moduli 
space is smooth in the (strong) sense of \cite[2.2]{BK}.

The second approach (\cite{Hin},\cite{ciokap}), takes $\AN$ as the 
category of differential $\Z$-graded local Artinian $\K$-algebras with 
residue field $\K$ that are homotopy equivalent to algebras 
concentrated in nonpositive degrees and the generalized Schlessinger 
condition are the above 1), 2), 3) together with the following: 
\begin{description}
\item[4)] $F$ sends quasiisomorphisms in $\AN$ into invertible maps 
in $\cat{Set}$.
\end{description}
This fourth condition is the principal responsible of some hidden 
smoothness  
phenomena of extended deformation functors. This approach works 
usually quite 
well for the extension of (classical) prorepresentable deformation 
functors, e.g. the Hilbert functor, but fails to be good in more 
general situations (see the discussion in  \cite[3.3]{K}).

In this paper we propose a notion of extended\footnote{From now on we 
always omit the term `extended' before `deformation functor'}
deformation functor which contains the previous two as 
specializations. In doing that we always keep in mind two 
general ideas: the first is taken from \cite{K}, while the second 
is nowadays standard.
\begin{description}
\item[a)] $\AN$ is the category of all differential $\Z$-graded 
associative (graded)-commutative 
local Artinian $\K$-algebras with residue field $\K$.
\item[b)] Over a field of characteristic 0, every reasonable 
deformation problem is governed by a differential graded Lie algebra; 
therefore the functor of Maurer-Cartan 
solutions of a differential graded Lie algebra modulo gauge action 
must be considered as the basic example of deformation functor.
\end{description}

Some easy examples show that the above condition 4) 
is not compatible with a) and b); 
this is one of the motivation of this paper and 
forces us to give a  more technical definition of 
deformation functor, see \ref{defofun}. 
The paper goes as follows:

In Section 1 we recall some definitions from rational homotopy theory  
and we fix the notation for the rest of the paper.

In Section 2 we propose a precise definition of deformation functors 
in terms of some Schlessinger's type conditions and we show that our 
notion is compatible with the above ideas a) and b). 
As in the classical case, every deformation functor $F$ has a  
tangent space $TF[1]$ which is a $\Z$-graded vector space.
For notational 
convenience we always replace local Artinian rings with their maximal 
ideals; this means that, in our description, $\AN$ is the category of finite 
dimensional associative commutative $\Z$-graded nilpotent differential 
$\K$-algebras. 

In Section 3 we show that the inverse function theorem holds for 
morphisms of deformation functors (Cor. \ref{CorIMT}); 
this implies in particular a 
conceptually easier proof of the well known fact that every 
quasiisomorphism of differential graded Lie algebras induces an 
isomorphisms of deformation functors. We also prove that every 
deformation functor has a natural complete obstruction theory with 
obstruction space $TF[2]$ (by definition $TF[2]$ is a graded vector 
space with $TF[2]^{i}=TF[1]^{i+1}$).

In Section 4 we generalize the main result of \cite{Sch}; 
every deformation functor $F$ induces naturally a functor $[F]$
defined in the homotopy category  $K(\AN)$ (see 
\ref{defhtp}). Theorem \ref{equivalencethm} asserts that, if the 
tangent space $TF[1]$, is finite dimensional then $[F]$ is 
prorepresented  by a pronilpotent differential  algebra $(R,d)$, where 
$R$ is the inverse limit of 
$R/R^{n+1}=\oplus_{i=1}^n(\odot^i TF[1]\dual)$ 
and $d$ is a differential of degree +1.

In Section 5 we will deal with extended 
deformation functors in 
the set-up of dg-coalgebras and we introduce the deformation functors 
of a $L_{\oo}$-algebra.

In Section 6 we will
prove the existence of a natural graded Lie algebra structure over the 
(shifted) tangent space of a deformation functor. 

In Section 7 we give a necessary and sufficient condition, always 
satisfied in concrete examples, for a deformation functor to be 
isomorphic to the deformation functor of a  $L_{\oo}$-algebra.

This paper owes its existence to 
the participation of the author to the ``Deformation quantization 
seminar'', held in Scuola Normale Superiore di Pisa during the academic 
year 1998-99. It is a pleasure to thank here E. Arbarello, G. Bini, A. 
Canonaco, P. de Bartolomeis, F. de Vita, D. Fiorenza, G. Gaiffi, M. 
Grassi, M. 
Polito and R. Silvotti for useful and stimulating discussions.

\subsection*{Notation}

We will always work over a fixed field $\K$ of characteristic 0.
All vector spaces, linear maps, algebras, tensor products, 
derivations etc.. are understood of being over $\K$, unless otherwise 
specified.

By $\Sigma_{n}$ we denote the symmetric group of permutations of
$\{1,\ldots,n\}$. By $S(p,q)\subset\Sigma_{p+q}$ we denote the set of 
unshuffles of type $(p,q)$; by definition $\sigma\in S(p,q)$ if and 
only if $\sigma(i)<\sigma(i+1)$ for every $i\not=p$. The unshuffles are 
a set of representative for the  cosets of 
$\Sigma_{p}\times \Sigma_{q}$ inside $\Sigma_{p+q}$.\\

\section{Differential graded algebras and homotopy}

We denote by: 
$\cat{Set}$ the category of sets in a fixed universe; 
$\cat{G}$  the category of $\Z$-graded vector 
spaces. If $V$ is a graded 
vector space and $v\in V$, $v\not=0$, is a homogeneous element, we denote 
by $\bar{v}\in \Z$ its degree.\\  
$\cat{DG}$  denotes the category of complexes of 
vector spaces, also called differential graded vector spaces.
Every  object in $\cat{DG}$ can be considered as a 
pair\footnote{For notational convenience, when the risk of confusion is 
remote, we prefer to write $V$ for the complex omitting the 
differential; the same will be made for all algebraic structures 
builded over a graded vector space.}
$(V,d)$, 
where $V$ is
an object in $\cat{G}$ and the differential  
$d\colon V\to V$ is a linear map such that $d(V^{i})\subset V^{i+1}$ 
and $d^2=0$.  
We will consider $\cat{G}$ as the full subcategory of $\cat{DG}$  whose 
objects are the complexes with zero differential.\\

For every complex of vector spaces $V\in \cat{DG}$ we denote, as usual, 
by $Z^{*}(V), B^{*}(V)$ and $H^{*}(V)$ the cocycles, coboundaries and 
cohomology of $V$. 

The tensor product in the category $\cat{DG}$  is 
defined in the standard way 
\[
(V\otimes W)^i=\oplus_{j\in\Z}V^j\otimes W^{i-j},\qquad d(v\otimes
w)=dv\otimes w+(-1)^{\bar{v}}v\otimes dw
\]
The twisting map $T\colon V\otimes W\to W\otimes V$ is defined by 
$T(x\otimes y)=(-1)^{\bar{x}\,\bar{y}}y\otimes x$ and then extended by 
linearity; $T$ is a morphism of complexes and define a $\Z/2\Z$-action
on $V\otimes V$. This action extends naturally to an action of the 
symmetric group $\Sigma_n$ on $\otimes^{n}V$.\\
The {\em Koszul} sign  
$\epsilon(\sigma;v_{1},\ldots,v_{n})=\pm 1$, $\sigma\in \Sigma_{n}$, 
is defined by
\[\sigma(v_{1}\otimes\ldots\otimes v_{n})=
\epsilon(\sigma;v_{1},\ldots,v_{n})
(v_{\sigma(1)}\otimes\ldots\otimes v_{\sigma(n)}).\]
In the next formulas we simply write $\epsilon(\sigma)$ instead of 
$\epsilon(\sigma;v_{1},\ldots,v_{n})$ when there is no ambiguity about
$v_1,\ldots, v_n\in V$.\\
The symmetric power $\odot^{n}V$ is the quotient of 
$\otimes^{n}V$ by the subspace generated by the elements
$v_{1}\otimes\ldots\otimes v_{n}-
\sigma(v_{1}\otimes\ldots\otimes v_{n})$.

Given an integer $n$, the shift\footnote{Here we follows the notation of 
\cite{K}; in \cite{Qui} the complex $V[-n]$ is called the `$n$-fold 
suspension of $V$'} 
functor 
$[n]\colon \cat{DG}\to \cat{DG}$ is defined by setting 
$V[n]=\K[n]\otimes V$, 
$f[n]=Id_{\K[n]}\otimes f$, 
where 
\[ 
\K[n]_i=\left\{\begin{array}{ll}
\K & \mbox{if $i+n=0$}\\
0 & \mbox{otherwise}\\
\end{array}\right.
\]

Given $v\in V$ we also denote by $v[n]=1_{\K[n]}\otimes v\in V[n]$, 
where $1_{\K[n]}$ is the unity of $\K$ shifted in degree $-n$; note 
that $\bar{v[n]}=\bar{v}-n$.

More informally, the complex $V[n]$ is the complex $V$ with indices 
shifted 
by $n$ and differential multiplied by $(-1)^n$. Note that 
$\Hom_{\cat{DG}}(V,W[n])$ are the linear maps $f\colon V\to W$ such that 
$f(V_i)\subset W_{i+n}$ for every $i$ and $f d_V=d_{W[n]} f=(-1)^n d_W f$. 
\par\smallskip

A commutative, associative graded algebra  is the data of a 
graded vector space $A\in \cat{G}$ together an associative,
linear multiplication map
$\colon A\otimes A\mapor{} A$ such that $A_iA_j\subset 
A_{i+j}$ and $ab=(-1)^{\bar{a}\,\bar{b}}ba$.

We denote by $\cat{GA}$ the category of (commutative, associative) 
graded algebras. A graded algebra $A$ is called nilpotent if $A^n=0$ for 
$n>>0$; clearly every nilpotent algebra is without unit.\\

A commutative associative differential graded algebra (dg-algebra for 
short)
is an object $(A,d)\in \cat{DG}$ such that $A$ is a graded algebras and 
$d$ is a derivation of degree 1. This means that $d$ satisfies the 
(graded) Leibnitz rule $d(ab)=d(a)b+(-1)^{\bar{a}}  a d(b)$. 
If $A$ has a unit 1, we assume moreover that $1\not\in d(A)$.

We denote by $\cat{DGA}$ the category of dg-algebras. 
In the sequel of the paper we consider 
$\cat{DG}$  
and $\cat{GA}$ as 
the full subcategories of $\cat{DGA}$ whose object are 
respectively the dg-algebras 
with trivial multiplication and with trivial differential. 
We also denote by $\cat{NDGA}$ the full 
subcategory of nilpotent dg-algebras and by $\AN$ the full subcategory of
nilpotent dg-algebras which are finite dimensional as $\K$-vector space.
Note that in our convention  $\cat{DG}\cap\AN$ denotes the dg-algebras 
$A\in\AN$ with trivial multiplication.

	\begin{ex}

If $A\in\AN$, then $\K\oplus A$ is a local Artinian ring; 
conversely it is easy to see that every local Artinian dg-algebra with 
unit and residue field $\K$ has the form $\K\oplus A$ for some $A\in \AN$
	
	\end{ex}
	
A module over a dg-algebra $A$ is the data of a complex $M\in \cat{DG}$ 
together two associative 
multiplication maps $A\otimes M\to M$ (left multiplication), 
$M\otimes A\to M$ (right multiplication) 
which are morphisms in $\cat{DG}$ commuting with the twisting map. This 
means that:
\begin{itemize}
\item 
$am=(-1)^{\bar{a}\,\bar{m}}ma$ for every homogeneous $a\in A$, 
$m\in M$.
\item $d_M(am)=d_A(a)m+(-1)^{\bar{a}}ad_M(m)$, 
~$d_M(ma)=d_M(m)a+(-1)^{\bar{m}}md_A(a)$.
\end{itemize}
If $A$ has a unit $1$ we also assume that $1\cdot m=m\cdot 1=m$.

If $M$ is an $A$-module and $n\in\Z$ 
then there is a natural structure of $A$-module 
over $M[n]$ having the same right multiplication and 
the left multiplication induced via the twisting map $T$.

A derivation of a dg-algebra $A$ into an $A$-module $M$ is a morphism of 
graded vector spaces $h\colon A\to M$ which satisfy the Leibnitz rule 
$h(ab)=h(a)b+ah(b)$. We denote by $\Der(A,M)$ the $A_0$-module 
of derivations $h\colon A\to M$.  Notice that the differential of $A$
is an element of $\Der(A,A[1])$ and that if $A^2=AM=0$ then
$\Der(A,M)$ is the space of morphisms in the category $\cat{G}$ from
$A$ to $M$.

In the  above set-up it is also defined an $A$-module $\Der^*(A,M)\in 
\cat{DG}$ with $\Der^n(A,M)=\Der(A,M[n])$, the obvious left multiplication 
and differential 
\[
\delta\colon \Der(A, M[n])\to\Der(A,M[n+1]);\qquad
\delta(h)=d_{M[n]}h-hd_A
\]

Given a morphism of $A$-modules $f\colon M\to A$  such that 
$f(M)M=0$ (in most applications 
$M$ will be a square-zero  ideal of the dg-algebra $A$ and $f$ the 
inclusion) we 
define the mapping cone as the dg-algebra $C=A\oplus M[1]$ with the 
product 
$(a,m)(b,n)=(ab,an+mb)$ (note that, as a graded algebra, $C$ is the 
trivial extension of $A$ by $M[1]$) and differential 
\[
d_C=\left(\begin{array}{cc}
d_A&f\\ 
0&d_{M[1]}\\
\end{array}\right)
\colon A\oplus M[1]\to A[1]\oplus M[2]
\]
The reader must be careful here: the product in $C$ is defined using the 
$A$-module structure of $M[1]$. We left to the reader the easy 
verification that the mapping cone $C$ is a dg-algebra, the inclusion 
$A\to C$ is a morphism of dg-algebras and the projection $C\to M[1]$ is a 
derivation.

Conversely, given a derivation $h\colon B\to N$ we define the 
derived inverse mapping cone as the dg-algebra $D=B\oplus N[-1]$ 
with product $(a,m)(b,n)=(ab,an+mb)$ 
and differential 
\[
d_D=\left(\begin{array}{cc}
d_B&0\\ h&d_{N[-1]}\\
\end{array}\right)
\colon B\oplus N[-1]\to B[1]\oplus N
\]

Here the projection $D\to B$ is a morphism of dg-algebras and the 
inclusion $N[-1]\to D$ is a morphism of $D$-modules.

We denote by $\K[t_1,\ldots,t_n,dt_1,\ldots,dt_n]$ the dg-algebra of polynomial 
differential forms on the affine space $\A^n$ with the de Rham 
differential. 
We have $\K[t,dt]=\K[t]\oplus\K[t]dt$ and 
\[
\K[t_1,\ldots,t_n,dt_1,\ldots,dt_n]=\tensor{i=1}{n}\K[t_i,dt_i]
\]
Since $\K$ has characteristic 0, it is immediate to see that 
$H_{*}(\K[t,dt])=\K[0]$ and then by K\"{u}nneth formula 
$H_*(\K[t_1,\ldots,t_n,dt_1,\ldots,dt_n])=\K[0]$.\\
Note that for every dg-algebra $A$ and every $s=(s_1,\ldots,s_n)\in\K^n$ 
we have an evaluation morphism $e_s\colon 
A\otimes\K[t_1,\ldots,t_n,dt_1,\ldots,dt_n]\to A$ defined by 
\[
e_s(a\otimes p(t_1,\ldots,t_n,dt_1,\ldots,dt_n))=p(s_1,\ldots,s_n,0,\ldots,0)a
\]

For every dg-algebra $A$ we denote $A[t,dt]=A\otimes\K[t,dt]$; if $A$ is 
nilpotent then $A[t,dt]$ is still nilpotent. In the next sections the 
obvious fact that $A[t,dt]$ does not belong to $\AN$ for every $A\not=0$ 
will cause some problems whose solution is the introduction, for every 
positive real number $\epsilon>0$ of the dg-subalgebra
\[
A[t,dt]_\epsilon=A\oplus \oplus_{n>0} A^{\lceil n\epsilon\rceil}\otimes
(\K t^n\oplus \K t^{n-1}dt)\subset A[t,dt]
\]
It is clear that if $A\in \AN$ then $A[t,dt]_\epsilon\in \AN$ for every 
$\epsilon>0$ and $A[t,dt]$ is the union of all $A[t,dt]_\epsilon 
$, 
$\epsilon>0$.

	\begin{defn}\label{defhtp}
Given two morphisms of dg-algebras $f,g\colon A\to B$, a homotopy
between $f$ and $g$ is a morphism $H\colon A\to B[t,dt]$ such 
that $H_0:=e_0\circ H=f$, $H_1:=e_1\circ H=g$ (cf. \cite[p. 
120]{GrMor}).\\ 
We denote by $[A,B]$ the quotient of $\Hom_{\cat{DGA}}(A,B)$ by the 
equivalence relation generated by homotopy. 
Since composition respects homotopy equivalence we can also consider 
the homotopy categories 
$K(\cat{DGA})$, $K(\cat{NDGA})$ and $K(\AN)$; it is convenient to think 
of $K(\cat{DGA})$ as the category with the same objects of $\cat{DGA}$ 
and with morphisms $Mor(A,B)=[A,B]$.
A morphism in 
$\cat{DGA}$ (resp.: $\cat{NDGA}$, $\AN$) is called a homotopy equivalence 
if it gives an isomorphism in the corresponding homotopy category. 

   \end{defn}

Two homotopic morphisms induce the same morphism in homology, see 
\cite[p. 120]{GrMor}.
If $A,B\in \cat{DG}$, then two morphisms $f,g\colon A\to B$ are 
homotopic in the sense of \ref{defhtp} if and only if $f$ is homotopic to 
$g$ in the usual sense. In particular every acyclic complex is contractible
as a dg-algebra.

\begin{defn}

An  extension in $\cat{NDGA}$ is a short exact sequence 
\begin{equation}
\label{smex}
0\mapor{}I\mapor{}A\mapor{\alpha}B\mapor{}0
\end{equation}
such that $\alpha$ is a morphism in $\cat{NDGA}$ and $I$ is an ideal of 
$A$ such that $I^{2}=0$. The extension (\ref{smex}) 
is called {\em small} if $AI=0$, it is called {\em acyclic small} 
if it is small and $I$ is an acyclic complex, or equivalently if
$\alpha$ is a quasiisomorphism.   				
\end{defn}
			
The easy proof of the following facts is left to the reader:
\begin{itemize}

\item  Every surjective morphism $A\mapor{\alpha} B$ in the category 
$\AN$ is the composition of a finite number of small extensions.   

\item  If $A\mapor{\alpha} B$ is a surjective quasiisomorphism in  
$\AN$ and $A_i=0$ for every $i>0$ then $\alpha$ is the 
composition of a finite number of acyclic small extensions. This is 
generally false if $A_{i}\not=0$ for some $i>0$. 

\item For every morphism $\alpha\colon A\to B$ in $\AN$ there exists a 
surjective homotopy equivalence $\gamma\colon C\to A$ in $\AN$ such
that $\alpha\gamma$ is homotopic to a surjective map.  If moreover
$A_{i}=B_{i}=0$ for every $i>0$ and $\alpha\colon H_{0}(A)\to
H_{0}(B)$ is surjective then it is possible to choose $C$ such that
$C_{i}=0$ for every $i>0$.

\item Let $F\colon K(\AN)\to \cat{Set}$ be a functor such that 
$F(\alpha)$ is bijective for every acyclic small extension $\alpha$; 
then for every $A,B\in \AN$ with 
$A_{i}=B_{i}=0$ for every $i>0$
and every quasiisomorphism $\gamma\colon A\to B$ the map $F(\gamma)$ 
is bijective.
\end{itemize}

\bigskip

\section{Extended deformation functors}	

  \begin{defn}\label{defofun}
A covariant functor $F\colon \AN \to 
\cat{Set}$ is called a {\em predeformation functor} if the following
conditions are satisfied:
\begin{enumerate}
\item $F(0)=\{0\}$ is the one-point set.
\item {\em (Generalized Schlessinger's conditions)}: For every pair of 
morphisms 
$\alpha\colon A\to C$, $\beta\colon B\to C$ in 
$\AN$ consider the map
\[\eta\colon F(A\times_C B)\to F(A)\times_{F(C)}F(B)\]
Then:
\begin{enumerate}
\item\label{schlessi1} $\eta$ is surjective when $\alpha$ is surjective.
\item\label{schlessi2} $\eta$ is bijective when $\alpha$ is 
surjective and $C\in \cat{DG}\cap\AN$ is acyclic.
\end{enumerate}
\item\label{quasismooth} {\em (Quasismoothness)}:
For every acyclic small extension
\[
0\mapor{} I\mapor{} A\mapor{}B\mapor{}0
\]
the induced map $\rho\colon F(A)\to F(B)$ is surjective.
\end{enumerate}
A predeformation functor $F$ is called a {\em deformation functor} if the 
map 
$\rho$ defined in \ref{quasismooth} is bijective.
	\end{defn}

The predeformation functors (resp.: deformation functors) together their 
natural transformations form a category which we denote by 
$\cat{PreDef}$ (resp.: $\cat{Def}$). Note that definition \ref{defofun}
also makes sense for covariant functors $F\colon\cat{NDGA}\to \cat{Set}$.

   \begin{lem}\label{eqschl}
For a covariant functor $F\colon\AN\to \cat{Set}$ with $F(0)=\{0\}$ it 
is sufficient to check condition \ref{schlessi2} of Definition 
\ref{defofun}  
when $C=0$ and when $B=0$ separately.
\end{lem}

\begin{pf} 

Follows immediately from the equality 
\[ A\times_CB=(A\times B)\times_C 0\]
where $A\mapor{\alpha}C$, $B\mapor{\beta}C$ are as in 
\ref{schlessi2} of \ref{defofun} and 
the fibred product on the right comes from the morphism $A\times 
B\to C$, $(a,b)\to \alpha(a)-\beta(b)$.\qed
\end{pf}

\begin{prop}\label{algefree}
Let $V$ be a graded vector space and let $(R,d)$ be a differential 
algebra such 
that, as a graded algebra, $R$ is the inverse limit of 
$R/R^{n+1}=\oplus_{i=1}^{n}(\odot^{i}V)$. 
Then the functor 
$h_{R}\colon\cat{NDGA}\to \cat{Set}$,  
\[h_{R}(A)=\Hom_{\cat{DGA}}(R/R^{n},A),\qquad n>>0\] 
is a predeformation functor.
\end{prop}

\begin{pf}
We have $F(0)=\{0\}$ and $F(A\times_CB)=F(A)\times_{F(C)}F(B)$; the 
only nontrivial condition to check is the surjectivity of 
$F(A)\to F(B)$ for every acyclic small extension $0\to I\to 
A\mapor{p}B\to 0$. For simplicity we prove this in the particular case 
$d(R)\subset R^2$; the proof in the general case is essentially the same 
but more messy.

Let $\phi\colon R/R^{n}\to B$ be a morphism of dg-algebras, $\{v_i\}\subset V$ 
a homogeneous basis and denote $b_i=\phi(v_i)$. For every set of 
liftings $a_i\in A$, $p(a_i)=b_i$, $\bar{a_i}=\bar{b_i}$, there exists 
a unique  morphism of graded algebras $\psi\colon R/R^{n+1}\to A$ such that 
$\psi(v_i)=a_i$; we need to prove that it is possible to choose the 
liftings 
such that $\psi(dv_i)=da_i$ for every $i$.

We first note that, since $AI=0$ and $\psi(dv_i)-da_i\in I$, 
the restriction $\psi\colon R^2/R^{n+1}\to A$ is 
a morphism of dg-algebras independent 
from the choice of the liftings $\{a_i\}$. In particular for every $i$, 
$d(\psi(dv_i)-da_i)=0$ and, being $I$ acyclic, there exist $s_i\in I$, 
$\bar{s_i}=\bar{a_i}$, such that $\psi(dv_i)-da_i=ds_i$. It is now 
sufficient to change $a_i$ with $a_i+s_i$ in order to transform $\psi$ 
into a morphism of differential graded algebras.\qed 
  \end{pf}

\begin{defn}
\label{defoquasismooth}
We shall call the differential algebra $R$ introduced in 
\ref{algefree}
algebraically free (or quasismooth) complete differential algebra. 
We also say that it is complete 
semifree (or semismooth) if 
there exists a filtration $0=H^0\subset H^1\subset....\subset V$ of 
graded vector spaces such that $\cup_{i=0}^{\oo} H^i=V$ and 
$d(H^{i+1})\subset\liminv\oplus_{j=1}^{n}(\odot^{j}
H^i$).\\
We shall say in addition that $R$ is minimal if $d(R)\subset R^{2}$. 
Note that, as a graded vector space, $V=R/R^{2}$. 
\end{defn}
The terms complete in \ref{defoquasismooth} is because $R$ is complete 
for the $R$-adic topology.

\begin{exs}~\\
\begin{enumerate}
\item If $R$ is a complete semifree algebra, it is easy to see that 
for every surjective quasiisomorphism $A\to B$ of nilpotent dg-algebras 
the 
morphism  $h_{R}(A)\to h_{R}(B)$ is 
surjective. This property is generally false if $R$ is assumed 
algebraically free.    

\item Let $(R,d)$ be a complete quasismooth differential algebra, 
$V=R/R^{2}$,
then $R$ is
semifree if either $V_i=0$ for every $i>0$ or $R$ is minimal, $V_1=0$
and $V_j=0$ for every $j<0$.  In particular the maximal ideal of the
algebra of functions of a smooth formal pointed dg-scheme (see \cite{ciokap})
is a complete semifree algebra.
\end{enumerate}
  \end{exs}

\begin{lem}\label{critdef}
For a predeformation functor $F\colon \AN\to \cat{Set}$ the following 
conditions are equivalent:
\begin{description}
\item[i)] $F$ is a deformation functor.
\item[ii)] $F$ induces a functor $[F]\colon K(\AN)\to \cat{Set}$.
\item[iii)] If $I\in \cat{DG}\cap \AN$ is acyclic then $F(I)=\{0\}$.
\end{description}
	\end{lem}
			
\begin{pf}
i)$\solose$ ii) Let $H\colon A\to B[t,dt]$ be a homotopy, we need to
prove that $H_0$ and $H_1$ induce the same morphism from $F(A)$ to
$F(B)$.  Since $A$ is finite-dimensional there exists $\epsilon>0$
sufficiently small such that $H\colon A\to B[t,dt]_{\epsilon}$; now
the evaluation map $e_0\colon B[t,dt]_{\epsilon}\to B$ is a finite
composition of acyclic small extensions and then, since $F$ is a
deformation functor $F( B[t,dt]_{\epsilon})=F(B)$; For every $a\in
F(A)$ we have $H(a)=iH_0(a)$, where $i\colon B\to B[t,dt]_{\epsilon}$
is the inclusion and then $H_1(a)=e_1H(a)=e_1iH_0(a)=H_0(a)$.\\
ii)$\solose$ iii) In the homotopy category every acyclic complex is 
isomorphic to 0.\\
iii)$\solose$ i) We need to prove that 
for every acyclic small extension
\[
0\mapor{} I\mapor{} A\mapor{\rho}B\mapor{}0
\]
the map $F(A)\to F(B)$ is injective, this is  done by showing that 
the diagonal map $F(A)\to F(A)\times_{F(B)}F(A)$ is surjective; in order 
to prove this it is sufficient to prove that the diagonal map 
$A\to A\times_BA$ induces a surjective map $F(A)\to F(A\times_BA)$. 
We have a canonical isomorphism $\theta\colon A\times I\to A\times_BA$, 
$\theta(a,x)=(a,a+x)$ which sends $A\times \{0\}$ onto the diagonal; since 
$F(A\times I)=F(A)\times F(I)=F(A)$ the proof is concluded.\qed
\end{pf}

A standard argument in Schlessinger's theory \cite[2.10]{Sch} 
shows that for every 
predeformation functor $F$ and every $A\in \AN\cap \cat{DG}$ there 
exists a natural structure of vector space on $F(A)$, where the sum and 
the scalar multiplication are described by the maps
\[
A\times A\mapor{+}A\quad \solose\quad F(A\times A)=F(A)\times 
F(A)\mapor{+}F(A)
\]
\[ 
s\in\K,\quad A\mapor{\cdot s}A\qquad \solose\qquad F(A)\mapor{\cdot s}F(A)
\]
For every
deformation functor $F$ and every integer $i$ we denote
\[ 
T^iF=F(\K[i-1]),\qquad TF=\oplus_{i\in\Z}T^iF[-i]
\]

Every natural transformation $\phi\colon F\to G$ of deformation
functors induces linear maps $T^iF\to T^iG$ and then a morphism of 
graded vector spaces $TF\to TG$.

\begin{defn}
Given a deformation functor $F$, the graded vector spaces $TF[1]$ 
and $TF[2]$ are called respectively  tangent and obstruction space 
of $F$.
\end{defn}

\begin{thm}\label{standardconstruction}

Let $F$ be a predeformation functor, then there exists a deformation 
functor $F^+$ and a natural transformation $\eta\colon F\to F^+$ such 
that for every deformation functor $G$ and every natural transformation 
$\phi\colon F\to G$ there exists unique a natural transformation 
$\psi\colon F^+\to G$ such that $\phi=\psi\eta$.
\end{thm}
\begin{pf}
We first define a functorial relation $\sim$ on the sets $F(A)$, $A\in 
\AN$; we set $a\sim b$ if and only if there exists $\epsilon>0$ and 
$x\in F(A[t,dt]_\epsilon)$ such that $e_0(x)=a$, $e_1(x)=b$. By 
\ref{critdef} 
if $F$ is a deformation functor then $a\sim b$ if and only if $a=b$. 
Therefore if we define $F^+$ as the quotient of $F$ by the equivalence 
relation generated by $\sim$ and $\eta$ as the natural projection, 
then there exists a unique $\psi$ as in the 
statement of the theorem. We only need to prove that $F^+$ is a 
deformation 
functor.  

\begin{step}{1} If $C\in\cat{DG}\cap\AN$ is acyclic  then $F^+(C)=\{0\}$.\\
Since $C$ is acyclic there exists a homotopy $H\colon C\to 
C[t,dt]_\epsilon$, $\epsilon\le 1$, 
such that $H_0=0$, $H_1=Id$; it is then clear that 
for every $x\in F(C)$ we have $x=H_1(x)\sim H_0(x)=0$.
\end{step}

\begin{step}{2} $\sim$ is an equivalence relation on $F(A)$ for every 
$A\in \AN$.\\
This is essentially standard  (see e.g. \cite[p. 125]{GrMor}). 
In view of the inclusion $A\to A[t,dt]_\epsilon$ the relation $\sim$ is 
reflexive. The symmetry is proved simply by remarking that the 
automorphism of dg-algebras
\[ A[t,dt]\to A[t,dt];\qquad a\otimes p(t,dt)\to a\otimes p(1-t,-dt)\]
preserves the subalgebras $A[t,dt]_\epsilon$ for every $\epsilon>0$.

Consider now $\epsilon>0$ and $x\in F(A[t,dt]_\epsilon)$, $y\in 
F(A[s,ds]_\epsilon)$ such that $e_0(x)=e_0(y)$; we need to prove that 
$e_1(x)\sim e_1(y)$.

Write $\K[t,s,dt,ds]=\oplus_{n\ge 0}S^n$, where $S^n$ is the $n$-th 
symmetric power of the acyclic complex $\K t\oplus\K s\mapor{d}\K 
dt\oplus\K ds$ and define $A[t,s,dt,ds]_\epsilon=A\oplus 
\oplus_{n>0} A^{\lceil n\epsilon\rceil}\otimes S^n$. There exists a 
commutative diagram

\[\begin{array}{ccc}
A[t,s,dt,ds]_\epsilon&\mapor{t=0}&A[s,ds]_\epsilon\\
\mapver{s=0}&&\mapver{s=0}\\
A[t,dt]_\epsilon&\mapor{t=0}&A\\
\end{array}\]

The kernel of the surjective morphism
\[ A[t,s,dt,ds]_\epsilon \mapor{\eta} A[t,dt]_\epsilon\times_A 
A[t,dt]_\epsilon\]
is equal to $\oplus_{n>0}A^{\lceil n\epsilon\rceil}\otimes 
(S^n\cap I)$, where $I\subset \K[t,s,dt,ds]$ is the homogeneous 
differential ideal generated by $st,sdt,tds, dtds$. Since $I\cap S^n$ is 
acyclic for every $n>0$, the morphism $\eta$ is a finite 
composition of acyclic small extensions.

Let $\xi\in F(A[t,s,dt,ds]_\epsilon)$ be a lifting of $(x,y)$ and let 
$z\in F(A[u,du]_\epsilon)$ be the image of $\xi$ under the morphism
\[ A[t,s,dt,ds]_\epsilon\to A[u,du]_\epsilon,\qquad 
t\to 1-u,\quad s\to u\]
The evaluation of $z$ gives $e_0(z)=e_1(x)$, $e_1(z)=e_1(y)$.
\end{step}

\begin{step}{3}
If $\alpha\colon A\to B$ is surjective then 
\[ F(A[t,dt]_\epsilon)
\mapor{(e_0,\alpha)}F(A)\times_{F(B)}F(B[t,dt]_\epsilon)\] 
is surjective.\\
		
It is not restrictive to assume $\alpha$ a small extension with kernel 
$I$. The kernel of $(e_0,\alpha)$ is equal to 
$\oplus_{n>0}(A^{\lceil n\epsilon\rceil}\cap I)\otimes (\K t^n\oplus \K 
t^{n-1}dt)$ and therefore  $(e_0,\alpha)$ is an acyclic small extension.
\end{step}

\begin{step}{4}
The functor $F^+$ satisfies \ref{schlessi1} of \ref{defofun}.\\

Let $a\in F(A)$, $b\in F(B)$ be such that $\alpha(a)\sim\beta(b)$; by 
Step 3 there exists $a'\sim a$, $a'\in F(A)$ such that 
$\alpha(a')=\beta(b)$ and then the pair $(a',b)$ lifts to $F(A\times_C 
B)$. 
\end{step}

\begin{step}{5}
The functor $F^+$ satisfies \ref{schlessi2} of \ref{defofun}.\\

By  \ref{eqschl} it is sufficient to verify the condition separately for 
the cases $C=0$ and $B=0$. 
When $C=0$ the situation is easy: in fact $(A\times B)[t,dt]_\epsilon=
A[t,dt]_\epsilon\times B[t,dt]_\epsilon$, 
$F((A\times B)[t,dt]_\epsilon)=
F(A[t,dt]_\epsilon)\times F(B[t,dt]_\epsilon)$ and the relation $\sim$
over $F(A\times B)$ is the product of the relations $\sim$ over $F(A)$ 
and $F(B)$; this implies that $F^+(A\times B)=F^+(A)\times F^+(B)$.

Assume now $B=0$, then the fibred product $D:=A\times_C B$ is equal to 
the kernel of $\alpha$. We need to prove that the map $F^+(D)\to F^+(A)$ 
is injective.
Let $a_0,a_1\in F(D)\subset F(A)$ and let $x\in F(A[t,dt]_\epsilon)$ be 
an element such that $e_i(x)=a_i$, $i=0,1$. Denote by $\bar{x}\in 
F(C[t,dt]_\epsilon)$ the image of $x$ by $\alpha$.
 
Since $C$ is acyclic there exists a morphism of graded vector spaces
$\sigma\colon C\to C[-1]$ such that $d\sigma+\sigma d=Id$ and we can 
define a morphism of complexes
\[ h\colon C\to (\K s\oplus\K ds)\otimes C\subset C[s,ds]_1;\qquad 
h(v)=s\otimes v+ds\otimes\sigma(v)\]
The morphism $h$ extends in a natural way to a morphism 
\[ h\colon C[t,dt]_\epsilon\to (\K s\oplus\K ds)\otimes C[t,dt]_\epsilon\]
such that for every scalar $\zeta\in\K$  there exists a commutative diagram

\[\begin{array}{ccc}
C[t,dt]_\epsilon&\mapor{h}&(\K s\oplus\K ds)\otimes C[t,dt]_\epsilon\\
\mapver{e_{\zeta}}&&\mapver{Id\otimes e_{\zeta}}\\
C&\mapor{h}&(\K s\oplus\K ds)\otimes C
\end{array}\]

Setting $\bar{z}=h(\bar{x})$ we have $\bar{z}_{|s=1}=\bar{x}$, 
$\bar{z}_{|s=0}=\bar{z}_{|t=0}=\bar{z}_{|t=1}=0$. By step 3 $\bar{z}$ 
lifts to an element $z\in F(A[t,dt]_\epsilon[s,ds]_1)$ such that 
$z_{|s=1}=x$; Now the specializations ${z}_{|s=0}$, ${z}_{|t=0}$, 
${z}_{|t=1}$ are annihilated by $\alpha$ and therefore give a chain of 
equivalences in $F(D)$
\[ a_0=z_{|s=1, t=0}\sim z_{|s=0, t=0}\sim z_{|s=0, t=1}\sim z_{|s=1, 
t=1}=a_1\]
proving that $a_0\sim a_1$ inside $F(D)$.
\end{step}

The combination of 
Steps 1, 4 and 5 tell us that $F^+$ is a deformation functor.\qed 
\end{pf}

For later use we point out that in the previous proof (Steps 2,3) 
we have also showed the following

\begin{lem}\label{smoothpredef}
Let $F\colon \AN\to\cat{Set}$ be a predeformation functor and let 
$\nu\colon F\to F^{+}$ be the projection onto the associated 
deformation functor.\\
Then for every surjective morphism $A\to B$ in $\AN$, the map 
\[ F(A)\to F(B)\times_{F^{+}(B)}F^{+}(A)\]
is surjective.
\end{lem}

For computational purposes it is often useful the following 

\begin{lem}\label{calcolotangente}
Let $S$ be a complex of vector spaces and assume that the functor 
\[ \cat{DG}\cap\AN\to\cat{Set};\ C\to Z^{1}(S\otimes C)\]
is the restriction of a predeformation functor $F$. Then for every 
complex $C\in\cat{DG}\cap\AN$ we have $F^+(C)=H^{1}(S\otimes C)$; in particular 
$T^iF^+=H^{i}(S)$.
\end{lem}

\begin{pf}
Let $C\in \cat{DG}\cap \AN$ be a fixed complex and set $A=S\otimes C$; 
since $A^2=C^2=0$ for every $\epsilon>0$ we have $A[t,dt]_\epsilon=S\otimes 
C[t,dt]_\epsilon$.

For every integer $i\ge 2$ let $p_i(t)\in\K[t]$ be a fixed monic 
polynomial of degree $i$ such that $p(0)=p(1)=0$; 
every element $\xi\in 
Z^{1}(A[t,dt]_\epsilon)=F(C[t,dt]_\epsilon)$ can be written as 
\[\xi=a_0+a_1t+\sum_{i\ge 2}a_ip_i(t)+b_1dt+\sum_{i\ge 
2}b_ip_i'(t)dt\]
with $a_i\in A_1$, $b_i\in A_0$, $da_0=0$ and $db_i=a_i$ for every $i$. 
In particular $e_1(\xi)-e_0(\xi)=a_1=db_1\in B^{1}(A)$.
Conversely if $a_0,a_1\in Z^{1}(A)$ and $a_1-a_0=db_1$ for some $b_1\in 
A_0$ then $a_j=e_j(a_0+(a_1-a_0)t+b_1dt)$ for $j=0,1$.
\qed
\end{pf}

\begin{ex}\label{2.11} 
Let $R$ be an algebraically free complete differential algebra;
\begin{enumerate}
\item	 For every 
complex $C$ we have 
\[h_{R}(C)=\Hom_{\cat{DG}}(R/R^{2},C)=
Z^{0}(\Hom^{*}(R/R^{2},\K)\otimes{C})\]
and then 
$T^{i}h_{R}^{+}=H^{i-1}(\Hom^{*}(R/R^{2},\K))$.

\item
The functor  
$[R,-]\colon\cat{NDGA}\to \cat{Set}$ is a deformation functor with 
tangent space 
\[T^{i}[R,-]=H^{i-1}(\Hom^{*}(R/R^{2},\K))
=H^{i-1}((R/R^{2})\dual).\]
The proof 
of this fact is essentially the same (but more easy) of 
\ref{standardconstruction} and it is 
omitted. If the graded vector space $R/R^{2}$ is 
finite-dimensional then it is easy to prove directly that 
$h_{R}^{+}=[R,-]$
 
\item The projection $h_{R}\to [R,-]$ factors through 
a morphism of deformation functors 
$h_{R}^{+}\to [R,-]$ which is an 
isomorphism on tangent spaces. 
\end{enumerate}

\end{ex}

All the deformation functors are quasismooth by definition; a stronger
notion of smoothness is given by the natural generalization of the
classical notion of smooth morphism.

\begin{defn}\label{defosmooth}
A natural transformation of deformation functors $\theta\colon F\to G$ 
is called smooth if, for every morphism $A\to B$ in $\AN$, which is 
surjective in homology  the natural map 
\[ F(A)\to F(B)\times_{G(B)}G(A)\]
is surjective. A deformation functor $F$ is called smooth if the 
trivial morphism $F\to 0$ is smooth.
\end{defn}

\begin{prop}
Let $(S,d)$ be a minimal complete quasismooth differential 
algebra. Then 
the deformation functor $[S,-]$ is smooth if and only if $d=0$.
\end{prop}

\begin{pf} If $d=0$ it is easy to see that $[S,-]$ is smooth.
Conversely assume $d\not=0$; let $n\ge 2$ be the greatest integer 
such that $d(S)\subset S^{n}$ and
let $\pi\colon S/S^{2}\to V$ be a finite dimensional quotient such that the
composition  $H\mapor{d}S^{n}/S^{n+1}\to \odot^{n}(V)$ is nonzero.

Let $(R,0)$ be the complete free differential algebra 
with $R/R^{2}=V$; 
we shall see that $[S,R/R^{n+1}]\to [S,R/R^{n}]$ is not 
surjective, this implies that $[S,-]$ is not smooth.
Since $S/S^{n}$ and $R/R^{n}$ have differential $d=0$ we have  
$[S, R/R^{n}]=\Hom_{\cat{DGA}}(S,R/R^{n})$ and therefore it is sufficient 
to show that the projection $\pi\colon S\to R/R^{n}$ does not lift 
to $R/R^{n+1}$. Assume  that there exists $f\colon S\to R/R^{n+1}$ 
which is a 
lifting of $\pi$; let $x\in S/S^{2}$ be such that $\pi(dx)\not=0$ in $R^{n}$, 
then
$f(x)=\pi(x)+\alpha$ for some $\alpha\in R^{n}/R^{n+1}$; in particular
the restriction of $f$ to $S^{2}$ is equal to $\pi$.  We must
therefore have $0=df(x)=f(dx)=\pi(dx)\not=0$ giving a
contradiction.\qed
\end{pf}

The second basic example of deformation functor is the deformation 
functor of a differential graded Lie algebra; we recall the well known 

\begin{defn}
A differential graded Lie algebra (DGLA in short terms) 
is a triple $(L,d,[,])$ with $(L,d)\in 
\cat{DG}$ and $[,]\colon L\times L\to L$ is a bilinear map such that:
\begin{itemize}
\item $[L^i, L^j]\subset L^{i+j}$.
\item $[x,y]+(-1)^{\bar{x}\,\bar{y}}[y,x]=0$.
\item (Jacoby identity) 
$[x,[y,z]]=[[x,y],z]+(-1)^{\bar{x}\,\bar{y}}[y,[x,z]]$
\item $d[x,y]=[dx,y]+(-1)^{\bar{x}}[x,dy]$
\end{itemize}
\end{defn}

\medskip

\begin{ex} Let $A$ be a dg-algebra, then the $A$-module $\Der^*(A,A)$ has a 
natural structure of DGLA with bracket 
$[\delta,\tau]=\delta\circ\tau-(-1)^{\bar{\delta}\bar{\tau}}\tau\circ\delta$, 
where the composition $\circ$ is made by considering $\delta$ and $\tau$ 
as linear endomorphism of the $\K$-vector space $A$.
\end{ex}

We shall denote by $\cat{DGLA}$ the category of differential graded Lie 
algebras: a morphism of DGLA is simply a morphism of complexes which 
commutes with brackets. 

The Maurer-Cartan elements of a DGLA $L$ are by definition 
\[ MC(L)=
\left\{x\in L^1\,\left|\, dx+\frac{1}{2}[x,x]=0\right.\,\right\}\]
Clearly every morphism of differential graded Lie algebras $L\to N$
sends the Maurer-Cartan elements of $L$ into the ones of $N$.

Given a DGLA $L$ and $A\in \AN$, the tensor product $L\otimes{A}$ has 
a 
natural structure of nilpotent DGLA with 
\[ (L\otimes{A})^i=\oplus_{j\in\Z}L^j\otimes A_{i-j}\]
\[ d(x\otimes a)=dx\otimes a+(-1)^{\bar{x}}x\otimes da\]
\[ [x\otimes a, y\otimes b]=(-1)^{\bar{a}\,\bar{y}}[x,y]\otimes ab\]
Every morphism of DGLA, $L\to N$ and every morphism $A\to B$ in $\AN$ give 
a natural commutative diagram of morphisms of differential graded Lie 
algebras
\[\begin{array}{ccc}
L\otimes{A}&\mapor{}&N\otimes{A}\\
\mapver{}&&\mapver{}\\
L\otimes{B}&\mapor{}&N\otimes{B}\end{array}\]

\begin{defn}\label{MCDGLA}
The Maurer-Cartan functor $MC_L\colon\AN\to\cat{Set}$ of a DGLA $L$ is 
given by 
\[ MC_L(A)=MC(L\otimes{A})\]
\end{defn}

\begin{lem}
In the notation above $MC_L$ is a predeformation functor; the 
resulting functor 
$MC\colon\cat{DGLA}\to \cat{PreDef}$ is faithful.
\end{lem}

\begin{pf}

It is evident that $MC_L(0)=0$ and for every pair of morphisms 
$\alpha\colon A\to C$, $\beta\colon B\to C$ in $\AN$ we have
\[ MC_L(A\times_C B)=MC_L(A)\times_{MC_L(C)}MC_L(B)\]
Let $0\mapor{} I\mapor{} A\mapor{}B\mapor{}0$ be an acyclic small 
extension and $x\in MC_L(B)$. Since $\alpha$ is surjective there 
exists $y\in (L\otimes{A})^1$ such that $\alpha(y)=x$. Setting 
\[ h=dy+\frac{1}{2}[y,y]\in (L\otimes{I})^2\]
we have 
\[dh=\frac{1}{2}d[y,y]=[dy,y]=[h,y]-\frac{1}{2}[[y,y],y].\] 
By 
Jacoby identity $[[y,y],y]=0$ and, since $AI=0$ also $[h,y]=0$; thus 
$dh=0$ and, being $L\otimes{I}$ acyclic by K\"{u}nneth formula, there 
exists 
$s\in (L\otimes{I})^1$ such that $ds=h$. The element $y-s$ lifts $x$ 
and 
satisfies Maurer-Cartan equation. We have thus proved that $MC_L$ is a 
predeformation functor.

The natural identification $L^i=MC_L((\K t\oplus\K dt)[i])$ shows that 
$MC\colon\cat{DGLA}\to \cat{PreDef}$ is faithful.\qed
\end{pf}

\begin{rem}\label{2.15}
It is easy to prove that a differential graded Lie algebra can be
recovered, up to isomorphism, from its Maurer-Cartan functor.
\end{rem}

In order to motivate our definition of deformation functor we point out 
that, if $A\to B$ is a surjective quasiisomorphism in $\AN$, then in 
general $MC_L(A)\to MC_L(B)$ is not surjective.
As an example take $\K$  algebraically closed and 
$L$ a finite-dimensional non-nilpotent  Lie 
algebra, considered as a DGLA concentrated in degree 0. Fix $a\in L$ 
such that $ad(a)\colon L\to L$ has an eigenvalue $\lambda\not=0$. Up to 
multiplication of $a$ by $-\lambda^{-1}$ we can assume $\lambda=-1$.
Let $V\subset L$ be the image of $ad(a)$, the linear map $Id+ad(a)\colon 
V\to V$ is not surjective and then there exists $b\in L$ such that the 
equation $x+[a,x]+[a,b]=0$ has no solution in $L$.

Let $u,v,w$ be indeterminates of degree 1 and consider the dg-algebras 
\[ B=\K u\oplus \K v,\qquad B^2=0,\, d=0\]
\[ A=\K u\oplus \K v\oplus\K w\oplus\K dw,\qquad uv=uw=dw,\, vw=0\]
The projection $A\to B$ is a quasiisomorphism but the element $a\otimes 
u+b\otimes v\in MC_L(B)$ cannot lifted to $MC_L(A)$. In fact, if there 
exists 
$\xi=a\otimes u+b\otimes v+x\otimes w\in MC_L(A)$, then 
\[ 0=d\xi+\frac{1}{2}[\xi,\xi]=(x+[a,x]+[a,b])\otimes dw\]
in contradiction with the previous choice of $a,b$.

\medskip

To every differential graded Lie algebra $(L,d,[,])$ we can
associate a new DGLA $(L_d, d_d, [,]_d)$ by setting $L_d=L\oplus\K d$,
where $d$ is considered as an indeterminate of degree 1, $d_d(a+\alpha
d)=d(a)$ and \[ [a+\alpha d, b+\beta d]_d=[a,b]+\alpha
d(b)-(-1)^{\bar{a}}\beta d(a)\] Note that $a\in MC(L)$ if and only if
$a\in L^1$ and $[a+d,a+d]_d=0$.

If $L^0=L^0_d$ is a nilpotent Lie algebra then the adjoint action 
$L^0\times L^1_d\mapor{[,]} L^1\subset L^1_d$ gives a group action 
$exp(L^0)\times L^1_d\to L^1_d$
preserving the affine hyperplane $L^1+\alpha d$ for every fixed 
$\alpha\in \K$.
As a consequence of Jacoby identity (see 
\cite{ManettiDGLA}, \cite{GoldMil1} for details) this action 
preserves the quadratic cone $\{x\in L^1_d\,|\, [x,x]_d=0\}$ and then 
$exp(L^0)$ acts, via the above identification, on the set $MC(L)$.

The restriction of the  adjoint action to $MC(L)$ is usually called 
``gauge action''. Note that if, for $a\in L^0$ and $b\in MC(L)$, we have
$[a,b+d]_d=[a,b]-da=0$ then $e^a b=b$; in particular th set 
$K_b=\{[b,c]+dc\,|\, c\in L^{-1}\}$ is a Lie subalgebra of $L^0$ and 
$exp(K_b)$ is contained in the stabilizer of $b$.

For every DGLA $L$ and every $A\in \AN$ we define $\Def_L(A)$ as the 
quotient of $MC(L\otimes{A})$ by the gauge action of the group 
$exp((L\otimes{A})^0)$. The gauge action commutes with 
morphisms in $\AN$ and with morphisms of differential graded Lie algebras 
and then the above definition gives a functor $\Def_L\colon\AN\to 
\cat{Set}$ called the {\em deformation functors of $L$}: in fact we have

\begin{thm}\label{tangenteperDGLA}
For every differential graded Lie algebra $L$,  
$\Def_L$ is a deformation functor with $T^i\Def_L=H^i(L)$.
\end{thm}

\begin{pf} 
If $C\in \cat{DG}\cap\AN$  then $L\otimes{C}$ is an abelian DGLA, 
$MC_{L}(C)=Z^1(L\otimes{C})$ and the gauge action is given by 
$e^a b=b-da$, $b\in (L\otimes{C})^1$, 
$a\in (L\otimes{C})^0$. This proves 
that $\Def_L(C)=H^1(L\otimes{C})$ and in particular 
$T^i\Def_L=H^1(L\otimes\K[i-1])=H^i(L)$.

Since $\Def_L$ is the quotient of a predeformation functor and 
$\Def_{L}(I)=0$ for every acyclic  $I\in \cat{DG}\cap\AN$,  
it is sufficient to 
verify the generalized Schlessinger's conditions. Let $\alpha\colon A\to 
C$, $\beta\colon B\to C$ morphism in $\AN$ with $\alpha$ surjective. 
Assume there are given $a\in MC_L(A)$, $b\in MC_L(B)$ such that 
$\alpha(a)$ and $\beta(b)$ give the same element in $\Def_L(C)$; then there 
exists $u\in (L\otimes{C})^0$ such that 
$\beta(b)=e^u\alpha(a)$. Let $v\in 
(L\otimes{A})^0$ be a lifting of $u$, 
changing if necessary $a$ with its gauge 
equivalent element $e^va$, we may suppose $\alpha(a)=\beta(b)$, the pair 
$(a,b)$ lifts to $MC_L(A\times_C B)$; this proves that the map 
\[ \Def_L(A\times_C B)\to \Def_L(A)\times_{\Def_L(C)}\Def_L(B)\]
is surjective.

If $C=0$ then the gauge action $exp((L\otimes(A\times B))^0)\times 
MC_L(A\times B)\to MC_L(A\times B)$ is the direct product of the gauge 
actions 
$exp((L\otimes{A})^0)\times 
MC_L(A)\to MC_L(A)$ and 
$exp((L\otimes{B})^0)\times 
MC_L(B)\to MC_L(B)$ and therefore $\Def_L(A\times B)=\Def_L(A)\times 
\Def_L(B)$.

Finally assume $B=0$, $C$ acyclic with $C^{2}=0$ and set $D=\ker\alpha\simeq 
A\times_C B$. Let $a_1,a_2\in MC_L(D)$, $u\in (L\otimes{A})^0$ be 
such 
that $a_2=e^u a_1$; we need to prove that there exists $v\in (L\otimes 
{D})^0$ such that $a_2=e^v a_1$.

Since $\alpha(a_1)=\alpha(a_2)=0$ and $L\otimes{C}$ is an abelian DGLA we 
have $0=e^{\alpha(u)}0=0-d\alpha(u)$ and then there exists $h\in 
(L\otimes{A})^{-1}$ such that $d\alpha(h)=-\alpha(u)$, $u+dh\in 
(L\otimes{D})^0$. 
If  $w=[a_1,h]+dh$ then $e^w a_1=a_1$, $e^u e^w a_1=e^v a_1=a_2$ 
where $v=u*w$ is determined by Baker-Campbell-Hausdorff formula.\\
We claim 
that $v\in L\otimes{D}$; in fact $v=u*w\equiv u+w\equiv 
u+dh\pmod{[L\otimes{A}, L\otimes{A}]}$ and since $A^2\subset D$ we have 
$v=u*w\equiv u+dh\equiv 0 \pmod{L\otimes{D}}$.\qed
\end{pf}

\begin{cor}\label{2.17}
For every differential graded Lie algebra $L$ the natural projection 
$MC_L\to 
\Def_L$ induces (by \ref{standardconstruction}) 
a morphism $MC_L^+\to \Def_L$ which is an isomorphism on tangent spaces.
\end{cor}

\bigskip

\section{Obstruction theory and the inverse function theorem for 
deformation functors}

\begin{thm}\label{IMT}
A morphism of predeformation functors $\theta\colon F\to G$ is an 
isomorphism if and only 
if $\theta\colon F(A)\to G(A)$ is 
a bijection for every complex $A\in \AN\cap \cat{DG}$.

\end{thm}

\begin{pf}

The proof  uses the natural generalization to the 
differential graded case of some standard techniques in Schlessinger's 
theory, cf. \cite{FM1}.

Let $\theta\colon F\to G$ be a fixed natural transformation of 
predeformation functors such that $\theta\colon F(A)\to G(A)$ is 
bijective for every  $A\in\cat{DG}\cap \AN$. 

\begin{step}{1}
For every small extension
\[ 0\mapor{} I\mapor{} A\mapor{\alpha} B\mapor{}0\]
and every $b\in F(B)$ we have either 
$\alpha^{-1}(b)=\emptyset$ or $\theta(\alpha^{-1}(b))=
\alpha^{-1}(\theta(b))$.

As in the classical case \cite{Sch}, 
there exists an isomorphism of dg-algebras
\[
A\times I\mapor{} A\times_B A;\qquad (a,t)\to (a,a+t)
\]
and then for every predeformation functor $E$ there exists a natural 
surjective map
\[
\vartheta_E: E(A)\times E(I)=E(A\times I)\to E(A)\times_{E(B)}E(A)
\]
This implies in particular that there exists a natural transitive action 
of the vector space $E(I)$ on the fibres of the map $E(A)\to E(B)$. 
Moreover this action commutes with natural transformations of functors.

In our case we have a commutative diagram
\[
\begin{array}{ccc}
F(A)&\mapor{\alpha}&F(B)\\
\mapver{\theta}&&\mapver{\theta}\\
G(A)&\mapor{\alpha}&G(B)\\
\end{array}
\]
and compatible transitive actions of the vector 
space $F(I)=G(I)$ on the fibres of 
the horizontal maps. This proves Step 1.
\end{step}

\begin{step}{2}
Let $A\mapor{\alpha}V$ be a surjective morphism in $\AN$ with 
$V\in\cat{DG}\cap\AN$ 
acyclic. Denoting by  
$\iota\colon B=\ker\alpha\to A$ the inclusion, 
then an element  $b\in G(B)$ lifts to 
$F(B)$ if and only if $\iota(b)$ lifts to $F(A)$.

In fact, by condition \ref{schlessi2} of \ref{defofun}, the inclusion 
$\iota$ gives bijection $F(B)=\alpha^{-1}(0)\subset F(A)$, 
$G(B)=\alpha^{-1}(0)\subset G(A)$. Since $F(V)=G(V)$ we have 
\[
\iota\theta(F(B))=\theta(F(A))\cap \alpha^{-1}(0)=\theta(F(A))\cap 
\iota(G(B))
\]  
\end{step}

\begin{step}{3} Let 
\[ 0\mapor{} I\mapor{\iota} A\mapor{\alpha} B\mapor{}0\]
be a small extension and $b\in F(B)$. Then $b$ lifts to $F(A)$ if and 
only if $\theta(b)$ lifts to $G(A)$.

The only if part is trivial, let's prove the if part. Let $\tilde{a}\in 
G(A)$ be a lifting of $\theta(b)$ and let $C=A\oplus I[1]$ be the mapping 
cone of $\iota$. The projection $\pi\colon C\to B$ is a small acyclic 
extension and then there exists $b_1\in F(C)$ which lifts $b$. By 
Step 1 we can assume that $\theta(b_1)=j(\tilde{a})$, where $j\colon A\to 
C$ is the inclusion.

If $ob\colon C\to I[1]$ denotes the projection then 
$ob(\theta(b_1))=0$ and since $F(I[1])=G(I[1])$ we also have $ob(b_1)=0$; 
by generalized Schlessinger's conditions $b_1$ lifts to $F(A)$.
\end{step}

\begin{step}{4}
For every $A\in\AN$ the map $\theta: F(A)\to G(A)$ is surjective.

By induction over $\dim_{\K}A$ we can assume that $F(B)$ goes 
surjectively onto $G(B)$ for every proper quotient $B=A/I$. Let 
\[ 0\mapor{} I\mapor{} A\mapor{\alpha} B\mapor{}0\]
be a small extension with $I\neq 0$, $\tilde{a}\in G(A)$ a fixed element 
and $b\in F(B)$ such that $\theta(b)=\alpha(\tilde{a})$. By Step 3 
$\alpha^{-1}(b)$ is not empty and then by Step 1 $\tilde{a}\in 
\theta(F(A))$.
\end{step}

\begin{step}{5}
Let $a\in F(A)$, for every surjective morphism $f\colon A\to B$ 
in the category $\AN$ we define 
\[ S_F(a,f)=\{\xi\in F(A\times_B A)\, | 
\, \xi\to (a,a)\in F(A)\times_{F(B)}F(A)\} \]
By definition, if $f$ is a small extension and $I=\ker f$ then $S_F(a,f)$ 
is naturally isomorphic to the stabilizer of $a$ under the action of 
$F(I)$ on the fibre $f^{-1}(f(a))$. It is also clear that 
$\theta(S_F(a,f))\subset S_G(\theta(a),f)$.
\end{step}

\begin{step}{6}
For every $a\in F(A)$ and every surjective morphism $f\colon A\to 
B$ the map $\theta\colon S_F(a,f)\to S_G(\theta(a),f)$ is surjective.

This is trivially true if $B=0$, 
we prove the general assertion by induction on $\dim_{\K}B$.
Let 
\[
0\mapor{}I\mapor{}B\mapor{\alpha}C\mapor{}0
\]
be a small extension with $I\not=0$, set  $g=\alpha f$ and denote by 
$h\colon A\times_C A\to I$ the surjective morphism of dg-algebras defined 
by $h(a_1,a_2)=f(a_1)-f(a_2)$; the kernel of $h$ is $A\times_B A$, let's 
denote by $\iota\colon A\times_B A\to A\times_C A$ the natural inclusion.

By generalized Schlessinger's conditions the maps

\[ F(A\times_B A)\to \ker(F(A\times_C A)\mapor{h}F(I));\qquad
S_F(a,f)\to S_F(a,g)\cap \ker h
\]
are surjective.
Let $\tilde{\xi}\in S_G(\theta(a),f)$ and let $\eta\in S_F(a,g)$ such 
that $\theta(\eta)=\iota(\tilde{\xi})$. Since $F(I)=G(I)$ we have 
$h(\eta)=0$ and then $\eta$ lifts to some $\xi_1\in S_F(a,f)$.
Let $D=(A\times_C A)\oplus I[-1]$ be the derived inverse mapping cone of 
$h$, it is immediate to check that the projection maps
\[
\pi_1\colon D\to A\times_C A,\qquad (h,\pi_2)\colon D\to I\oplus I[-1]
\]
are surjective morphisms of dg-algebras, $I\oplus I[-1]$ is an acyclic 
complex and the kernel of $(h,\pi_2)$ is exactly $A\times_B A$; again by 
generalized 
Schlessinger's conditions $F(A\times_B A)\subset F(D)$, $G(A\times_B 
A)\subset G(D)$. Moreover there exists a cartesian diagram
\[
\begin{array}{ccc}
(A\times_B A)\times I[-1]&\mapor{}&A\times_B A\\
\mapver{}&&\mapver{\iota}\\
D&\mapor{\pi_1}&A\times_C A\\
\end{array}
\]
which gives surjective maps
\[
F(A\times_B A)\times F(I[-1]) \mapor{\varrho}F(A\times_B 
A)\times_{F(A\times_C A)}F(D)
\]
\[
G(A\times_B A)\times G(I[-1]) \mapor{\varrho}G(A\times_B 
A)\times_{G(A\times_C A)}G(D)
\]
This implies the existence of 
$v\in G(I[-1])=F(I[-1])$ such that 
$\varrho(\theta(\xi_1),v)=(\theta(\xi_1),\tilde{\xi})$; defining  
$\xi\in F(D)$ by the formula 
$\varrho(\xi_1,v)=(\xi_1,\xi)$ we get  $\theta(\xi)=\tilde{\xi}$ and 
then, by Step 2, $\xi\in S_F(a,f)$.
\end{step}

\begin{step}{7} 
For every $A\in \AN$ the map $\theta\colon F(A)\to G(A)$ is 
injective. 

This is true by assumption if $A^2=0$; if $A^2\neq 0$ we can suppose by 
induction that there exists a small extension 
\[ 0\mapor{} I\mapor{\iota} A\mapor{\alpha} B\mapor{}0\]
with $I\neq 0$ and $\theta\colon F(B)\to G(B)$ injective.

Let $a_1,a_2\in F(A)$ be two elements such that 
$\theta(a_1)=\theta(a_2)$; by assumption $f(a_1)=f(a_2)$ and then there 
exists $t\in F(I)$ such that  $\vartheta_F(a_1,t)=(a_1,a_2)$, 
$\vartheta_G(\theta(a_1),\theta(t))=(\theta(a_1),\theta(a_2))$ and then 
$\theta(t)\in S_G(\theta(a_1),\alpha)$. By Step 6 there exists $s\in 
S_F(a_1,\alpha)$ such that $\theta(s)=\theta(t)$ and by injectivity of 
$\theta\colon F(I)\to G(I)$ we get $s=t$ and then $a_1=a_2$.
\end{step}
~~~\qed

\end{pf}

The result of Theorem~\ref{IMT} is particularly useful for morphisms of 
deformation functors: In fact we have

\begin{lem}\label{Fpercomplessi}
Let $F\colon \AN\to \cat{Set}$ be a deformation functor; for every 
$I\in\cat{DG}\cap\AN$ there exists a natural isomorphism 
\[F(I)=\somdir{i\in\Z}{}\,TF[1]^i\otimes H_{-i}(I)=
\somdir{i\in\Z}{}T^{i+1}F\otimes H_{-i}(I)=
H^{0}(\Hom^{*}(I\dual,TF[1])).\]
\end{lem}

\begin{pf}
Let $s\colon H_*(I)\to Z_*(I)$ be a linear section of the natural 
projection, then the composition of $s$ with the natural embedding 
$Z_*(I)\mapor{\iota} I$ is unique up to homotopy and its cokernel is an 
acyclic complex, therefore it  
gives a well defined isomorphism  $F(H_*(I))\to F(I)$. This says that it 
is not restrictive to prove the lemma for complexes with zero differential.
Moreover since $F$ commutes with direct sum of complexes we can reduce to 
consider the case when $I\simeq \K^s[n]$ is a vector space concentrated 
in degree $-n$.
Every $v\in I$ gives a morphism $F(\K[n])\mapor{v}F(I)$ and we can define 
a 
natural map $T^{1+n}F\otimes I\to F(I)$, $x\otimes v\to v(x)$. It is easy 
to verify that this map is an isomorphism of vector spaces.\qed
\end{pf}

As an immediate consequence we have:

\begin{cor}\label{CorIMT}
A morphism of deformation functors $\theta\colon F\to G$ is an 
isomorphism if and only if it gives an isomorphism of tangent spaces
$\theta\colon TF[1]\mapor{\simeq}TG[1]$
\end{cor} 

\begin{cor}
For every differential graded Lie algebra $L$ there exists a natural 
isomorphism  $MC_L^+=\Def_L$.
\end{cor}

\begin{pf}
According to \ref{standardconstruction} the projection $MC_{L}\to 
\Def_{L}$ induces a natural transformation $MC_{L}^{+}\to 
\Def_{L}$ which is an isomorphism on tangent spaces by \ref{2.17}.\qed 
\end{pf}

\begin{cor}
Every quasiisomorphism $\phi\colon 
L\to N$ of differential graded Lie algebras  
induces an isomorphism  $\Def_{L}\simeq \Def_{N}$ between deformation 
functors.
\end{cor}

\begin{pf}
Trivial consequence of  \ref{tangenteperDGLA} and \ref{CorIMT}.\qed
\end{pf}

\begin{cor}  For every 
complete quasismooth  dg-algebra $R$ there exists an isomorphism 
$h_{R}^{+}=[R,-]$.  
\end{cor}
\begin{pf}
Follows from \ref{2.11} and \ref{CorIMT}. Note that this is a 
nontrivial result 
when the graded vector space $R/R^{2}$ is infinite 
dimensional.\qed
\end{pf}

The argument used in the proof of \ref{IMT} can be used to show the 
existence of a complete natural (and probably universal, see \cite{FM1} 
for the definition) obstruction theory 
for every deformation functor $F$.

Given a small extension
\[
e:\quad 0\mapor{}I\mapor{\iota}A\mapor{\alpha}B\mapor{}0
\]
we can define an ``obstruction map'' $ob_e\colon F(B)\to F(I[1])
=H^{0}(TF[2]\otimes{I})$ in the following way:\\
Let $C=A\oplus I[1]$ be the mapping cone of the inclusion $\iota$; since 
the projection $C\to B$ is an acyclic small extension, the 
projection $C\to I[1]$ gives a map $ob_e\colon F(B)=F(C)\to F(I[1])$.

The obstruction maps satisfy the following properties:
\begin{itemize}
\item $ob_e(b)=0$ if and only if $b$ lifts to $F(A)$.
\item (naturality) The obstruction maps commute with natural 
transformation of functors.
\item (base change) Given a morphism of small extensions
\[
\begin{array}{cccccccccc}
\,e:&\quad 0&\mapor{}&I&\mapor{}&A&\mapor{\alpha}&B&\mapor{}&0\\
\mapver{}&&&\mapver{f}&&\mapver{}&&\mapver{f}&&\\
\,e':&\quad 0&\mapor{}&I'&\mapor{}&A'&\mapor{\alpha'}&B'&\mapor{}&0\\
\end{array}
\]
we have $f\circ ob_e=ob_{e'}\circ f\colon F(B)\to F(I'[1])$.
\end{itemize}
The last two items are straightforward, while the first 
follows by generalized Schlessinger's conditions applied to the cartesian 
diagram 
\[
\begin{array}{ccc}
A&\mapor{}&0\\
\mapver{}&&\mapver{}\\
C&\mapor{}&I[1]\\
\end{array}
\]

\begin{ex}
If $0\mapor{}I\mapor{}J\mapor{}L\mapor{}0$ is a short exact sequence of 
complexes then the obstruction map $F(L)\to F(I[1])$ is the tensor 
product of the identity over $TF$ and the connecting homomorphism 
$H_*(L)\to H_*(I[1])$.
\end{ex}

\begin{rem}
The above construction also shed light on the obstruction theory 
for general extensions. 
In this case it is necessary to replace obstruction maps 
with obstruction sections. If 
\[
e:\quad 0\mapor{}I\mapor{\iota}A\mapor{\alpha}B\mapor{}0;\qquad I^2=0
\]
is an extension, $C=A\oplus I[1]$ is the mapping cone of the inclusion 
and $B\oplus I[1]$ is the trivial extension (i.e. the 
mapping cone of the zero map $I[1]\mapor{0}B$)  
then the obstruction section $ob_e\colon F(B)=F(C)\to F(B\oplus I[1])$ is 
the map induced by the natural projection $C\to B\oplus I[1]$. The same 
proof as above shows that $b\in F(B)$ lifts to $F(A)$ if and only if 
$ob_e(b)=j(b)$, where $j\colon F(B)\to F(B\oplus I[1])$ is the ``zero 
section'' induced by the natural inclusion $B\subset B\oplus I[1]$. 
Clearly if $e$ is a small extension then $BI=0$, $F(B\oplus 
I[1])=F(B)\times 
F(I[1])$ and we recover the notion of obstruction maps.
\end{rem}

\begin{ex}\label{primobs} (Primary obstruction maps)\\
Let $u,v$ be indeterminates of degrees $-i,-j$ respectively and 
$F\colon\AN\to\cat{Set}$ a deformation functor; let 
$(u,v)\subset\K[u,v]$ be the ideal generated by $u,v$ 
and denote by 
\[Q^1_{ij}\colon T^{1+i}F\times T^{1+j}F=F(\K u\oplus K v)\to 
F(\K uv[1])=T^{2+i+j}F\] 
the obstruction map associated to the small extension 
\[ 0\mapor{}\K uv\mapor{}\frac{(u,v)}{(u^{2},v^2)}\mapor{}\K u\oplus K 
v\mapor{}0.\]
\end{ex}

In the realm of obstruction theory belong the following results 
which we will 
use in the  next sections.\\
Consider a generic small extension in $\AN$
\[
e:\quad 0\mapor{}I\mapor{\iota}A\mapor{\alpha}B\mapor{}0
\]
with obstruction map $ob_e\colon F(B)\to F(I[1])$. Let's denote by 
$d_0\colon A\to A[1]$ the differential of $A$. 

For every morphism of 
dg-algebras $\phi\colon B\to I[1]$ the map 
$d_\phi=d_0+\iota\phi\alpha\colon 
A\to A[1]$ is a new differential on $A$: in fact $BI=0$, this implies that
$\iota\phi\alpha$ is a derivation and  
\[ d_\phi^2=d_0\iota\phi\alpha+\iota\phi\alpha d_0=
\iota(d_{I}\phi+\phi d_B)\alpha=\iota(-d_{I[1]}\phi+\phi 
d_B)\alpha=0.\]
Denote by $A_\phi$ the graded algebra $A$ with the differential $d_\phi$
and by 
\[
e_\phi:\quad 0\mapor{}I\mapor{\iota}A_\phi\mapor{\alpha}B\mapor{}0
\]
the corresponding small extension.

\begin{lem}\label{torsorobstruction}
In  the notation above $ob_{e_\phi}=ob_e+\phi\colon F(B)\to F(I[1])$.
\end{lem}

\begin{pf}
Let $A_\phi\oplus I[1]$, $A\oplus I[1]$ be respectively the mapping cones 
of the inclusions $I\mapor{\iota} A_\phi$, $I\mapor{\iota}A$ and define 
the following isomorphism of graded algebras
\[
h=
\left(\begin{array}{cc} Id&0\\  \phi\alpha&Id\\  \end{array}\right)
\colon A_\phi\oplus I[1]\to A\oplus I[1]
\]
Since
\[
\left(\begin{array}{cc} Id&0\\  \phi\alpha&Id\\  \end{array}\right)
\left(\begin{array}{cc} d_\phi&\iota\\  0&d_{I[1]}\\  \end{array}\right)
=
\left(\begin{array}{cc} d_0&\iota\\  0&d_{I[1]}\\  \end{array}\right)
\left(\begin{array}{cc} Id&0\\  \phi\alpha&Id\\  \end{array}\right)
\]
the morphism $h$ is also an isomorphism of dg-algebras.
There exists a commutative diagram
\[
\begin{array}{ccc}
A_\phi\oplus I[1]&\mapor{h}&A\oplus I[1]\\
\mapver{\alpha\times Id}&&\mapver{\alpha\times Id}\\
B\times I[1]&\mapor{l}&B\times I[1]\\
\end{array}\quad\qquad
l=\left(\begin{array}{cc} Id&0\\  \phi&Id\\  \end{array}\right)
\]
showing that, via the bijections $\alpha\colon F(A_\phi\oplus 
I[1])\mapor{\simeq}F(B)$, $\alpha\colon F(A\oplus 
I[1])\mapor{\simeq}F(B)$ the isomorphism $h$ acts trivially on $F(B)$. 
Therefore the obstruction section 
$Id\times ob_{e_\phi}:F(B)\to F(B)\times F(I[1])$ 
is precisely the 
composition with $l$ of the obstruction section $Id\times ob_{e}$.\qed
\end{pf}

Let's consider two dg-algebras $(B,d_B)\in\AN$, $(I,d_I)\in\AN\cap 
\cat{DG}$ and a small extension of nilpotent graded algebras

\[
\quad 0\mapor{}I\mapor{\iota}A\mapor{\alpha}B\mapor{}0
\]
   
Assume that there exists a lifting of $d_B$ to a derivation $d\in Der(A,A[1])$ 
such that $d_{|I}=d_I$. Such a lifting is not unique, any two liftings 
differ by an derivation of the form $\iota\phi\alpha$ for some 
derivation $\phi\colon B\to I[1]$.

\begin{lem} In the above notation,  
the derivation $d^2=\frac{1}{2}[d,d]\colon A\to A[2]$ 
may be factored as $d^{2}=\delta\alpha$, where 
$\delta\colon B\to I[2]$ is a 
morphism of dg-algebras  whose homotopy class 
$[\delta]\in [B,I[2]]$ is independent from the choice of the lifting. 
Moreover 
$[\delta]=0$ if and only if there exists a lifting $d$ as above 
such that ${d}^2=0$.
\end{lem}

\begin{pf}
It is evident that $d^2(A)\subset I[2]$,  $d^2(I)=0$ and then $d^2$ 
induces a derivation $\delta\colon B\to I[2]$; since $d^2\circ 
d=d_I\circ d^2=d_{I[2]}\circ d^2$ we have that $\delta$ is also a 
morphism of dg-algebras and then gives a morphism of complexes 
$\bar{\delta}\colon B/B^2\to I[2]$.

If $d'$ is another lifting as above and $\phi=d'-d$, then the 
derivation $\phi$ comes in the obvious way 
from a morphism of graded vector spaces $\bar{\phi}\colon B/B^2\to 
I[1]$; since
$ \delta'-\delta=(d+\phi)^2-d^2=d\phi+ \phi d$ we have that
$\bar{\phi}$ is the required homotopy between 
$\bar{\delta}$ and $\bar{\delta'}$.
If $\delta\colon B\to I[2]$ is homotopic to 0 then, 
being $I[2][t,dt]^2=0$, also $\bar{\delta}$ is homotopic to 0 and then 
there exists a 
derivation $\phi\colon A\to B\to B/B^2\to I[1]$ such that $d^2=d\phi+ 
\phi d$, $(d-\phi)^2=0$.\qed
\end{pf}

In the above set-up, for every deformation functor $F$, the map 
$\delta\colon F(B)\to F(I[2])$ is, a priori, 
independent from the choice of the 
lifting of the differential. 
A posteriori we have the following stronger result:

\begin{prop}\label{vanishing}
In the above notation $\delta(F(B))=0$; in particular if there exists 
$\xi\in F(B)$ such that the map 
\[ [B,I[2]]\to F(I[2]);\qquad [f]\to f(\xi)\]
is injective then there exists a dg-algebra structure on $A$ making
$\quad 0\mapor{}I\mapor{\iota}A\mapor{\alpha}B\mapor{}0$ a small 
extension in $\AN$.
\end{prop}

\begin{pf}
We have already seen that $\delta$ induces a morphism of complexes 
$\bar{\delta}\colon B/B^2\to I[2]$.
It is sufficient to construct a commutative diagram in $\AN$
\[
\begin{array}{ccc}
C&\mapor{\gamma}&B\\
\mapver{}&&\mapver{}\\
V&\mapor{\bar{\gamma}}&B/B^2\\
\end{array}
\]
such that $\gamma$ is an acyclic small extension, 
$V^{2}=0$ and $\bar{\delta}\bar{\gamma}(H_*(V))=0\subset H_*(I[2])$.

Our solution is to define $C=A\times I[1]$ as a graded algebra and 
$\gamma$ as the composition of $\alpha$ with the projection on the first 
factor. Over $C$ we put the differential 
\[ d_C=\left(
\begin{array}{cc}
d&\iota\\ -d^2& d_{I[1]}\\
\end{array}\right)\colon A\times I[1]\to A[1]\times I[2]
\]

A simple verification shows that $d_C$ is a derivation, $d_C^2=0$ and 
$\gamma$ is a morphism of dg-algebras. The kernel of $\gamma$ is the 
mapping cone of the identity $I\to I$ and therefore $\gamma$ is an 
acyclic small extension.

Finally we set $V=B/B^2\times I[1]$ and the differential $d_V$ is forced 
by $d_C$ to be equal to 
\[ d_V=\left(
\begin{array}{cc}
d_{B/B^2}&0\\  -\bar{\delta}& d_{I[1]}\\
\end{array}\right)
\]

Since $\bar{\delta}\bar{\gamma}(b,t)=\bar{\delta}(b)$ and 
$d_V(b,t)=0\solose\bar{\delta}(b)=d_{I[1]}(t)=-d_{I[2]}(t)$ the morphism 
$\bar{\delta}\bar{\gamma}$ is trivial in homology.\qed

\end{pf}

\section{Finite deformation functors}
\label{sec:finitedef}

The notion of homotopy equivalence of morphisms works well for
finitely generated or nilpotent dg-algebras but it seems quite restrictive
for general differential algebras.  In the set-up of projective limits of
nilpotent dg-algebras it is useful to introduce a notion which arises 
naturally when we consider projective limits in homotopy categories. 

\begin{defn}\label{prohomot}
Let $f,g\colon S\to R$ be morphisms of differential algebras, with $R$ 
complete for the $R$-adic topology (i.e. $R=\liminv R/R^{n}$); we 
shall say that $f$ is prohomotopy equivalent to $g$ if 
the morphisms $f_{n}, g_{n}\colon S\to R\to R/R^{n}$, 
composition of $f,g$ with the natural projection $R\to R/R^{n}$,
are homotopy equivalent for every $n$.
\end{defn}

It is clear that the relation defined in \ref{prohomot} 
is an equivalence relation, if $R$ is nilpotent it is the same of the 
usual homotopy equivalence. 

If $f,g\colon S\to R$ are prohomotopy equivalent morphisms of 
quasismooth complete differential algebras then for every $A\in \AN$ and every 
$\phi\colon R\to A$ there exists a factorization $\phi\colon R\to 
R/R^{n}\to A$, $n>>0$, and therefore the composition $\phi f, \phi g$ 
are homotopy equivalent. This means that $f$ and $g$ induces the 
same morphism between deformation functors $[R,-]\to [S,-]$.

Another aspect of prohomotopy is the following 

\begin{thm}\label{modellominimale}
Every complete quasismooth differential algebra is prohomotopy equivalent 
to a minimal complete quasismooth differential algebra. Two minimal complete 
quasismooth algebras are prohomotopy equivalent if and only if 
they are isomorphic.
\end{thm}

\begin{pf} As in the nongraded case, a morphism $f\colon S\to R$ of 
complete quasismooth differential algebras is an isomorphism if and only if 
$f_{2}\colon S/S^{2}\to R/R^{2}$ is an isomorphism. Since the cohomology 
of $(S/S^{2})\dual$ is isomorphic, up to a shift, to the tangent space of 
the 
deformation functor $[S,-]$, if $f$ is a prohomotopy equivalence  then 
$[S,-]\simeq [R,-]$ and $f_{2}$ is a quasiisomorphism; this proves 
that every prohomotopy equivalence of minimal complete quasismooth 
differential algebras is an isomorphism.

Let $(R, d)$ be a fixed complete quasismooth differential algebra and let 
$d^{1}_{1}\colon R/R^{2}\to R/R^{2}$ be the differential induced by 
$d$. Let $R/R^{2}=H\oplus W$ be a decomposition with $W$ 
acyclic complex and $d^{1}_{1}(H)=0$; fix homogeneous basis 
$\{h_{j}\}$ of $H$ and $\{v_{i},w_{i}\}$ of $W$ such that 
$d^{1}_{1}(v_{i})=w_{i}$ for every $i$. Since 
$d(v_{i})=w_{i}+\phi_{i}$ for some $\phi_{i}\in R^{2}$, up to the 
analytic change of coordinates 
\[ h_{j}\to h_{j},\quad v_{i}\to v_{i},\quad w_{i}\to dv_{i}\]
it is not restrictive to assume $dW\subset W$. In particular the 
ideal $(W)\subset R$ is a differential ideal and the quotient 
$S=R/(W)$ is a minimal complete quasismooth differential algebra.

Our aim is to show that the projection $R\mapor{\pi}S$ is a 
prohomotopy equivalence; as a first step we prove that there exists a 
right inverse $\gamma\colon S\to R$ and then we prove that 
$\gamma\pi\colon R/R^{n}\to R/R^{n}$ is homotopic to the identity for 
every 
$n>0$.

Since $R$ is complete it is sufficient to find a sequence of morphisms  
$\gamma_{n}\colon S\to R/R^{n}$, $n\ge 2$, such that $\gamma_{2}\colon 
S\to H=S/S^{2}\to R/R^{2}=H\oplus W$ is the natural inclusion and the 
diagrams
\[\begin{array}{ccc}
S&\mapor{\gamma_{n}}&R/R^{n}\\
&\searrow&\mapver{\pi}\\
&&S/S^{n}\\
\end{array}\qquad\qquad
\begin{array}{ccc}
S&\mapor{\gamma_{n+1}}&R/R^{n+1}\\
\mapver{\gamma_{n}}&\swarrow&\\
R/R^{n}&&\\
\end{array}\]
are commutative.

We note that the natural morphism 
\[ R/R^{n+1}\to R/R^{n}\times_{S/S^{n}}S/S^{n+1}\]
is an acyclic small extension for every $n\ge 2$. According to 
\ref{algefree} we may define inductively $\gamma_{n+1}\colon S\to 
R/R^{n+1}$ as a lifting of 
\[(\gamma_{n},p_{n+1})\colon S\to R/R^{n}\times_{S/S^{n}}S/S^{n+1}\]
where $p_{n}\colon S\to S/S^{n}$ denotes the projection.

Up to a change of coordinates of the form $h_{j}\to \gamma\pi{h_{j}}$, 
$v_{i}\to v_{i}$, $w_{i}\to w_{i}$ we can assume 
that $S\subset R$
is the complete subalgebra generated by $H$.
The homotopies 
\[ H_{n}\colon R\to R/R^{n}[t,dt],\quad h_{j}\to h_{j},\quad v_{i}\to 
v_{i}\otimes t,\quad w_{i}=dv_{i}\to d(v_{i}\otimes t)\]
gives the required prohomotopy between $\gamma\pi$ and the 
identity.\qed 
\end{pf}

Let's denote by $K(\operatorname{quasismooth})$ the 
category of quasismooth complete dg-algebras with morphisms up to 
prohomotopy equivalence; then there exists the prorepresantability
functor 
\[(K(\operatorname{quasismooth}))^{opp}\to \cat{Def},\qquad 
S\to [S,-]\]

\begin{defn}\label{finitedefofun}
A deformation functor $F\colon \AN\to \cat{Set}$ is called {\em 
finite} if its tangent space $TF[1]$ has finite total dimension.
\end{defn}

\begin{defn} A complete quasismooth dg-algebra $R$ is called finite 
if $[R,-]$ is a finite deformation functor; i.e. if the graded vector 
space $H_{*}(R/R^{2})$ has finite total dimension.
\end{defn}

Note that a minimal complete quasismooth dg-algebra $R$ is finite if 
and only if $\K\oplus R$ is Noetherian.
If $\cat{Def^{\flat}}\subset \cat{Def}$ denotes the full subcategory of 
finite 
deformation functors and 
$K(\operatorname{quasismooth})^{\flat}\subset 
K(\operatorname{quasismooth})$ 
denotes the 
full subcategory of finite complete quasismooth dg-algebras we have

\begin{thm}\label{equivalencethm}
The restriction of the prorepresentability functor 
\[\vartheta\colon (K(\operatorname{quasismooth})^{\flat})^{opp}
\to\cat{Def^{\flat}}\] is an equivalence of categories. 
\end{thm}

\begin{pf} 
After theorem \ref{modellominimale} it is  not restrictive to 
assume that every complete quasismooth dg-algebra is minimal 
and finite.

\begin{step}{1} $\vartheta$ is injective on morphisms.\\
Let $f,g\colon S\to R$ be morphisms of dg-algebras inducing the same 
natural transformation of functors $[R,-]\to [S,-]$; if $\pi_{n}\colon 
R\to R/R^{n}$ denotes the natural projection then 
$\pi_{n}f,\pi_{n}g\colon S\to R/R^{n}$ are homotopic maps and then 
$f$ and $g$ are prohomotopic by definition.\end{step}

\begin{step}{2} $\vartheta$ is surjective on morphisms.\\
Let $S,R$ be finite complete minimal quasismooth dg-algebras
and $\alpha\colon [R,-]\to [S,-]$ be a natural transformation of 
functors. In the same notation of Step 1 we will construct 
recursively a coherent 
sequence of dg-algebra morphisms $\gamma_{n}\colon S\to R/R^{n}$ such 
that $[\gamma_{n}]=\alpha[\pi_{n}]\in [S,R/R^{n}]$. The resulting 
inverse limit $\gamma=\liminv \gamma_{n}\colon S\to R$ will be a 
morphism of dg-algebras inducing $\alpha$.

Assume $\gamma_{1},\ldots,\gamma_{n}$ as above are constructed for a 
fixed $n$ and let $f\colon S\to R/R^{n+1}$ be a representative of 
$\alpha[\pi_{n+1}]$. Denoting by $p\colon R/R^{n+1}\to R/R^{n}$ the 
projection there exists a homotopy $H\colon S\to R/R^{n}[t,dt]$ such 
that $H_{0}=pf$, $H_{1}=\gamma_{n}$. As $\Hom_{\cat{DGA}}(S,-)$ is a 
predeformation functor, the argument of Step 3 in the proof of 
\ref{standardconstruction} shows that there exists a 
homotopy $K\colon S\to R/R^{n+1}[t,dt]$ such that $pK=H$ and 
$K_{0}=f$; it is therefore sufficient define 
$\gamma_{n+1}=K_{1}$.\end{step} 

\begin{step}{3} $\vartheta$ is surjective on isomorphism classes.\\
Let $F$ be a deformation functor with finite dimensional tangent 
space $TF[1]$. By \ref{CorIMT} we need to prove that 
there exists a minimal complete 
quasismooth dg-algebra $(R,d)$ and a natural transformation of 
functors 
$\xi\colon [R,-]\to F$ which is an isomorphism on tangent spaces.

For every integer $i$, let $V_{-i}$ be the dual of the vector space 
$TF[1]^{i}$, then $V=\oplus_{i}V_{i}=\Hom^{*}(TF[1],\K)$ is a finite 
dimensional graded vector space. Thinking of $V$ as an 
object in $\cat{DG}\cap\AN$, we have 
by \ref{Fpercomplessi} a natural isomorphism
\[
F(V)=\somdir{i}{}(TF[1]^{i}\otimes V_{-i})=\Hom_{\cat{DG}}(TF[1],TF[1]).
\]
We denote by $\xi_{2}\in F(V)$ the elements corresponding to the 
identity on $TF[1]$. Again by \ref{Fpercomplessi} for every  
$I\in \cat{DG}\cap\AN$, the map
\[ \xi_{2}\colon [V,I]\to F(I),\qquad [f]\to f(\xi_{2})\]
is bijective.
Let $R$ be the 
inverse limit of $R/R^{n+1}=\oplus_{i=1}^{n}(\odot^{i}V)$;
we want to define two coherent sequences, 
the first  of square zero differentials $d_{n}\colon R/R^{n}\to 
R/R^{n}[1]$ with $d_{2}=0$ and the second of elements 
$\xi_{n}\in F(R/R^{n},d_{n})$
lifting $\xi_{2}$. This will give a structure of minimal complete 
quasismooth dg-algebra 
on $R=\lim R/R^{n}$ and the required natural transformation 
\[ [R,-]\to F,\qquad [f\colon R\to A]\to f(\xi_{n}),\, n>>0.\]
By induction assume there are defined $d_{n},\xi_{n}$, since $\xi_{n}$ 
lifts $\xi_{2}$, by \ref{vanishing} and \ref{torsorobstruction} there 
exists a unique square-zero differential $d_{n+1}\colon R/R^{n+1}\to 
R/R^{n+1}[1]$ lifting $d_{n}$ and such that the obstruction to 
lifting $\xi_{n}$ to $(R/R^{n+1},d_{n+1})$ vanishes.\qed 
\end{step}
\end{pf}

\begin{cor}\label{rapprecor}
For every finite deformation functor $F\colon\AN\to\cat{Set}$, the 
induced  functor $[F]\colon K(\AN)\to\cat{Set}$ is prorepresentable by a 
quasismooth minimal complete differential algebra $R$. If $T^{i}F=0$ for every 
$i\le 0$ then $\K\oplus R$ is the algebra of functions of a formal 
pointed noetherian semismooth dg-scheme.
\end{cor}
\begin{pf} Since prorepresentable means  representable by a projective 
limit
the proof comes immediately from Theorem~\ref{equivalencethm}.\qed    
\end{pf}

\bigskip

\section{Differential graded coalgebras and $L_{\oo}$-algebras}

\begin{defn}\label{coalgebra}
A cocommutative coassociative $\Z$-graded coalgebra is the data of 
a graded vector space $C=\oplus_{n\in\Z}C^{n}$ and of a 
coproduct $\Delta\colon C\to C\otimes C$ such that:
\begin{itemize}
\item $\Delta$ is a morphism of graded vector spaces.
\item (coassociativity) $(\Delta\otimes Id_{C})\Delta=
(Id_{C}\otimes\Delta)\Delta\colon C\to C\otimes C\otimes C$.
\item (cocommutativity) $T\Delta=\Delta$.
\end{itemize}
A morphism of coalgebras is a morphism of graded vector spaces that 
commutes with coproducts.
\end{defn}
By coassociativity we can define the iterated coproduct 
$\Delta^{n}=(Id_{C}\otimes \Delta^{n-1})\Delta\colon C\to 
\otimes^{n}C$. A coassociative coalgebra $(C,\Delta)$ is called 
{\em nilpotent} if $\Delta^{n}=0$ for $n>>0$.

\begin{ex}\label{symmetric-coalgebra}
The reduced symmetric coalgebra is by definition 
$\bar{S(V)}=\oplus_{n>0}\odot^{n}V$, with the coproduct:
\[\Delta(v_{1}\odot\ldots\odot v_{n})=
\sum_{r=1}^{n-1}\sum_{\sigma\in 
S(r,n-r)}\epsilon(\sigma) 
(v_{\sigma(1)}\odot\ldots\odot v_{\sigma(r)})\otimes
(v_{\sigma(r+1)}\odot\ldots\odot v_{\sigma(n)}).\]
An easy computation about unshuffles shows that the map 
$N\colon\bar{S(V)}\to\bar{T(V)}$ is a morphism of coalgebras.\\	
The coalgebra $\bar{S(V)}$ is coassociative, cocommutative without 
counity. Note also that $V=\ker\Delta$. 
We shall say that a subcoalgebra $C\subset \bar{S(V)}$ is {\em 
homogeneous} if $C=\oplus_n (C\cap \odot^n V)$. 
\end{ex}

Following \cite{K} we denote by $C(V)$ the coalgebra 
$(\bar{S(V[1])},\Delta)$.

\begin{defn}\label{coderivation}
Let $(C,\Delta)$ be a coalgebra. 
A morphism  $d\colon C\to C[n]$
of graded vector spaces is called a {\em coderivation of degree $n$} 
if 
\[\Delta d=(d\otimes Id_{C}+Id_{C}\otimes d)\Delta.\]
A coderivation $d$ is called a codifferential if $d^{2}=d\circ d=0$.\\
More generally if $\theta\colon C\to D$ is a morphism of coalgebras, 
a morphism of graded vector spaces $d\colon C\to D[n]$ is called a 
coderivation (with respect to $\theta$) if 
\[\Delta_{D} d=(d\otimes \theta+\theta\otimes d)\Delta_{C}.\]

\end{defn}
Note that, by the rule of sign, if $d\colon C\to C[n]$ is a morphism 
of graded vector spaces then $(Id_{C}\otimes d)(x\otimes y)=
(-1)^{n\bar{x}}x\otimes d(y)$.
The coderivations of degree $n$ of a coalgebra $C$ form a vector space 
denoted by $\Coder^{n}(C,C)$.

\begin{lem}\label{coderperdue}
Let $\theta_{1},\theta_{2}\colon C\to D$ be morphisms of coalgebras 
and assume that there exists a graded subspace $B\subset C$ such that 
$\Delta_{C}(C)\subset B\otimes B$ and $\theta_{1}(b)=\theta_{2}(b)$ 
for every $b\in B$.
Then every coderivation $C\to D[n]$ with respect to $\theta_{1}$ is 
also a coderivation with respect to $\theta_{2}$.
\end{lem}
\begin{pf} Evident.\qed\end{pf}

The reduced symmetric coalgebra is a free object in the the category 
of graded, locally nilpotent, cocommutative, coassociative 
coalgebras. In fact holds the following:

\begin{prop}\label{Universalpropsym}
Let $V$ be a graded vector space, $\pi\colon \bar{S(V)}\to V$ the 
projection and $C$ a cocommutative coalgebra which is union of 
nilpotent subcoalgebras.
\begin{enumerate}
\item The composition with $\pi$ gives a bijection
between the set of coalgebra morphisms $\theta\colon C\to\bar{S(V)}$ 
and morphisms of graded vector spaces $m\colon C\to V$. The inverse 
is given by the formula
\[ \theta=\sum_{n=1}^{\oo}\frac{1}{n!}m_{n}\Delta_{C}^{n-1},\]
where $m_{n}(c_{1}\otimes\ldots\otimes c_{n})=m(c_{1})\odot\ldots\odot 
m(c_{n})$.
\item For a fixed coalgebra morphism $\theta\colon C\to \bar{S(V)}$,
the composition with $\pi$ gives an isomorphism 
\[ \Coder^{n}(C,\bar{S(V)})\to \Hom_{\cat{G}}(C,V[n]).\]
\end{enumerate}\end{prop}
\begin{pf}  1) is proved in \cite[Prop. 4.2 on page 285]{Qui}.\\
2) follows easily by 1) in view of the following observation: 
consider the coalgebra $C\oplus C[-n]$ with coproduct 
\[ \Delta(c_{1}+c_{2}[-n])=\Delta(c_{1})+\Delta(c_{2})[-n]+
T(\Delta(c_{2})[-n]).\]
(In the above formula we used the natural isomorphism $(C\otimes 
C)[-n]=C[-n]\otimes C$). 
Then $d\colon C[-n]\to \bar{S(V)}$ is a coderivation (with respect to 
$\theta$) if and only if $\theta+d\colon C\oplus C[-n]\to\bar{S(V)}$ 
is a morphism of coalgebras.\qed\end{pf}

\begin{defn}\label{DGC} 
By a  dg-coalgebra we intend a triple $(C,\Delta,d)$, where 
$(C,\Delta)$ is a
graded coassociative cocommutative 
coalgebra  and $d\colon C\to C[1]$ is a codifferential of degree 1. 
The category of dg-coalgebras, where morphisms are morphisms of 
coalgebras commuting with codifferentials, is denoted by $\cat{DGC}$. 
\end{defn}

\begin{ex}\label{duale-algebra}
Given $A=\oplus A_{i}\in\AN$ 
let $C=A\dual$ be its graded dual: i.e. $C=\oplus C^{i}$, where 
$C^{i}=\Hom_{\K}(A_{-i},\K)$. We denote by $<,>\colon C\times A\to 
\K$ the induced pairing.\\
The pairing $<c_{1}\otimes c_{2},a_{1}\otimes 
a_{2}>=<c_{1},a_{1}><c_{2},a_{2}>$ gives a natural isomorphism
$C\otimes C=(A\otimes A)\dual$ commuting with the twisting maps $T$ 
and we may define 
$\Delta$ as the transpose of the multiplication map $\mu\colon 
A\otimes A\to A$.\\
Then $(C,\Delta)$ is a coassociative cocommutative coalgebra. Note 
that $C$ is nilpotent if and only if $A$ is nilpotent.\\
Defining the codifferential in $C$ as the transpose of the 
differential in $A$ we get a finite dimensional nilpotent dg-coalgebra.
\end{ex}
\begin{defn}\label{L-infinito}
Let $V$ be a graded vector space; a codifferential of degree 1 over 
the symmetric coalgebra $C(V)=\bar{S(V[1])}$ is called an {\em $L_{\oo}$ 
structure} on $V$. The dg-coalgebra $(C(V),Q)$ is called an 
{\em $L_{\oo}$-algebra}.\\
A $L_{\oo}$-algebra $(C(V),Q)$ is called {\em minimal} if $Q^{1}_{1}=0$. 
\end{defn}

Let $(C(V),Q)$
be a $L_{\oo}$-algebra and let $Q^{i}_{j}\colon 
\odot^{j}(V[1])\to \odot^{i}(V[1])[1]$ 
be the composition of the codifferential $Q$ with 
the inclusion $\odot^{j}(V[1])\subset C(V)$ and the projection 
$C(V)[1]\to \odot^{i}(V[1])[1]$; by Proposition~\ref{Universalpropsym} 
the codifferential $Q$ is uniquely determined by the linear maps 
$Q^{1}_{j}$, $j>0$ by the explicit formula (cf. \cite{penkava})
\[Q(v_{1}\odot\ldots\odot v_{n})=
\sum_{k=1}^{n}\sum_{\sigma\in S(k,n-k)}\epsilon(\sigma)
Q^{1}_{k}(v_{\sigma(1)}\odot\ldots\odot v_{\sigma(k)})\odot 
v_{\sigma(k+1)}\odot\ldots\odot v_{\sigma(n)}.\]
In particular we have 
$Q^{i}_{j}=0$ for every $i>j$ and then the subcoalgebras 
$\oplus_{i=1}^{r}(\odot^{i}V)$ are preserved by $Q$.

We define a functor 
$MC_{V}\colon\AN\to\cat{Set}$
by setting: 
\[ MC_{V}(A)=\Hom_{\cat{DGC}}(A\dual, C(V))\]

We use the following terminology: for every coalgebra $C$ and every 
linear map $\phi\colon C\to C(V)[n]$ we denote by $\phi^{i}\colon C\to 
\odot^{i}(V[1])[n]$ the composition of $\phi$ with the projection 
$C(V)[n]\to \odot^{i}(V[1])[n]$.

\begin{prop}\label{MCL-infinito}
For every $L_{\oo}$-algebra $(C(V),Q)$, 
the functor $MC_{V}$ is a predeformation functor.
\end{prop}
\begin{pf} It is evident that $MC_{V}$ is a covariant functor and
$MC_{V}(0)=0$. Let $\alpha\colon A\to C$, $\beta\colon B\to C$ be 
morphisms in $\AN$, then it is easy to see that 
the dg-coalgebra $(A\times_{C}B)\dual$ is the fibered coproduct, in the 
category $\cat{DGC}$, of the morphisms
\[\begin{array}{ccc}
C\dual&\mapor{{}^{t}\!\alpha}&A\dual\\
\mapver{{}^{t}\!\beta}&&\\
B\dual&&\end{array}\]
and then 
\[ MC_{V}(A\times_{C}B)=MC_{V}(A)\times_{MC_{V}(C)}MC_{V}(B).\]
Let $0\mapor{}I\mapor{}A\mapor{}B\mapor{}0$ be a small acyclic 
extension in $\AN$, we want to prove that $MC_{V}(A)\to MC_{V}(B)$ is 
surjective.\\
We have a dual exact sequence
\[0\mapor{}B\dual\mapor{}A\dual\mapor{}I\dual\mapor{}0,\qquad
B\dual=I^{\perp}.\]
Since $IA=0$ we have $\Delta_{A\dual}(A\dual)\subset B\dual\otimes 
B\dual$.\\
Let $\phi\in MC_{V}(B)$ be a fixed element and $\phi^{1}\colon B\dual\to 
V[1]$; by Proposition~\ref{Universalpropsym} $\phi$ is uniquely 
determined by $\phi^{1}$. Let $\psi^{1}\colon A\dual\to V[1]$ be an 
extension of $\phi^{1}$, then, by \ref{Universalpropsym}, 
$\psi^{1}$ is induced by a unique 
morphism of coalgebras $\psi\colon A\dual\to C(V)$.\\
The map $\psi d_{A\dual}-Q\psi\colon A\dual\to C(V)[1]$ is a 
coderivation and then, 
setting $h=(\psi d_{I\dual}-Q\psi)^{1}\in Hom_{\cat{G}}(I\dual,
V[2])$,  we have that $\psi$ is a 
morphism of dg-coalgebras if and only if $h=0$.\\
Note that $\psi^{1}$ is defined up to elements of 
$Hom_{\cat{G}}(I\dual, V[1])=
(V[1]\otimes I)^{0}$ and, since $\Delta_{A\dual}(A\dual)\subset B\dual\otimes 
B\dual$, $\psi^{i}$ depends only by $\phi$ for every $i>1$.
Since $I$ is acyclic and $h d_{I\dual}+Q^{1}_{1}h=0$ there exists 
$\xi\in Hom_{\cat{G}}(I\dual, V[1])$ such that 
$h=\xi d_{I\dual}-Q^{1}_{1}\xi$ and then $\theta^{1}=\psi^{1}-\xi$ induces a 
dg-coalgebra morphism $\theta\colon A\dual\to C(V)$ extending $\phi$.
\qed\end{pf}

\begin{defn}\label{MCLinfty}
Let $(C(V),Q)$ be a $L_{\oo}$-algebra and let $\Def_{V}=MC_{V}^{+}$ 
be the deformation functor associated to the predeformation functor 
$MC_{V}$ (Theorem~\ref{standardconstruction}). 
We shall call $\Def_{V}$ the 
{\em deformation functor associated to the $L_{\oo}$-algebra 
$(C(V),Q)$}.   
\end{defn}

A morphism of $L_{\oo}$-algebras $C(V)\to C(W)$ induces in the obvious 
way a natural transformation $MC_{V}\to MC_{W}$ and then, by 
Theorem~\ref{standardconstruction}, 
a morphism $\Def_{V}\to \Def_{W}$.

For every $A\in \AN$
there exists an algebra structure on 
$\Hom_{\cat{G}}(A\dual,C(V))$ induced by the
natural isomorphism 
$\Hom_{\cat{G}}(A\dual,C(V))=(C(V)\otimes A)^{0}$; the product of  
$f,g\colon A\dual\to C(V)$ is given by
$fg=\mu(f\otimes g)\Delta_{A\dual}$, where $\mu\colon C(V)\otimes 
C(V)\to C(V)$ is the multiplication of the reduced symmetric algebra
$C(V)=\bar{S(V[1])}$. According to 
\ref{Universalpropsym} a morphism $f\colon A\dual\to C(V)$ is a 
morphism of coalgebras if and only if there exist
$m\in \Hom_{\cat{G}}(A\dual,V[1])$ such that 
\[ f=\exp(m)-1\:=\sum_{n=1}^{\oo}\frac{1}{n!}m^{n}.\]

\begin{lem}\label{MCequazione}
Let $Q$ be a $L_{\oo}$ structure on $V$ and 
$m\in (V[1]\otimes A)^{0}$. Then $\exp(m)-1$ belongs to 
$MC_{V}(A)=\Hom_{\cat{DGC}}(A\dual,C(V))$, via the isomorphism 
$\Hom_{\cat{G}}(A\dual,C(V))=(C(V)\otimes A)^{0}$, if and 
only if 
\[Id_{V[1]}\otimes d_{A}(m)=Q^{1}\otimes Id_{A}(\exp(m)-1)=
\sum_{n=1}^{\oo}\frac{1}{n!}Q^{1}_{n}\otimes Id_{A}(m^{n}).\]
\end{lem}

\begin{pf} Denote $\theta=\exp(m)-1\colon A\dual\to C(V)$, then 
$\theta d_{A\dual}$, $Q\theta$ are coderivations and then 
$\theta d_{A\dual}=Q\theta$ if and only if 
\[Id_{V[1]}\otimes d_{A}(m)=(\theta d_{A\dual})^{1}=Q^{1}\theta=
Q^{1}\otimes Id_{A}(\exp(m)-1).\]
\qed\end{pf}

The above lemma allow to define 
$MC_{V}(A)$ 
for every nilpotent dg-algebra $A$; therefore we can define, as in 
\cite[4.5.2]{K}, 
the deformation functor $Def_V$ as the quotient of $MC_V$ by the 
equivalence relation given by: $x,y\in MC_{V}(A)$, $x\sim y$ if and only 
if there exists a solution $z\in MC_{V}(A[t,dt])$  
such that $x=z_{t=0}$, $y=z_{t=1}$.
The equivalence of this definition with $MC_V^+$ follows immediately from 
the equality $A[t,dt]=\cup_{\epsilon>0}V\otimes 
A[t,dt]_{\epsilon}$, which implies 
that $MC_{V}(A[t,dt])$ is the direct limit of 
$MC_{V}(A[t,dt]_{\epsilon})$,
and the explicit construction of $MC_V^+$ given in 
Theorem~\ref{standardconstruction}.

\begin{lem}\label{Ti} 
Let $(C(V),Q)$ be a $L_{\oo}$-algebra and consider $V$ as a 
differential graded vector space with differential $Q^1_1$.
Then for every  $A\in \AN\cap \cat{DG}$ we have
\begin{enumerate}
\item $MC_{V}(A)=\Hom_{\cat{DG}}(A\dual, V[1])=Z^{1}(V\otimes A)$,
\item  $\Def_{V}(A)=H^{1}(V\otimes A)$ and then  $T^{i}\Def_{V}=H^{i}(V)$.
\end{enumerate}
\end{lem}
\begin{pf} Since $A^{2}=0$ we have $\Delta_{A\dual}=0$ and then 
$\phi^{i}=0$ for every $i\ge 2$ and every morphism of coalgebras
$\phi\colon A\dual\to C(V)$: this proves 1).
The proof of 2) follows from 1) and Lemma~\ref{calcolotangente}.\qed\end{pf}

The following result is well known (cf. \cite{K}, \cite{LadaMarkl}):

\begin{prop} Given a graded vector space $V$,
there exists a bijection between the set of DGLA structures 
$(V,d, [,])$ and 
the set of $L_{\oo}$-algebra structures  
$(C(V),Q)$ such that $Q^{1}_{j}=0$ for every $j\ge 3$.\\ 
Explicitly the bijection is given by
$d(w)[1]=-Q^{1}_{1}(w[1])$ and 
$[w_{1},w_{2}]=(-1)^{\bar{w_{1}}}Q^{1}_{2}(
w_{1}[1]\odot  w_{2}[1])$ 
\end{prop}

An easy computation shows that for every DGLA $(V,d,[,])$ the 
definitions \ref{MCDGLA} and \ref{MCLinfty} give the same functor
$MC_{V}$, cf. \cite{K}.

\bigskip

\section{The Lie algebra structure on $TF$}

Let $F$ be a deformation functor, $T=TF$.
For every finite dimensional graded subspace $V\subset T[1]$ we have 
$F(V\dual)=\Hom_{\cat{G}}(V,T[1])$ and therefore there exists a 
distinguished element $\xi_{V}\in F(V\dual)$ corresponding to the 
inclusion $V\subset T[1]$.

\begin{prop} 
There exists a unique morphism of graded vector spaces 
$Q^{1}_{2}\colon \odot^{2}(T[1])\to T[2]$ with the following 
property:\\
For every pair of finite dimensional graded subspaces $V\subset 
H\subset T[1]$ such that $Q^{1}_{2}(\odot^{2}V)\subset H$ 
consider the subcoalgebra  $A=H\oplus(\odot^{2}V)\subset C(T)$ 
endowed 
with the codifferential
\[\left(\begin{array}{cc}
0&Q^{1}_{2}\\ 0&0\end{array}\right)\colon
H\somdir{}{}(\symm{}{2}V)\to (H\somdir{}{}(\symm{}{2}V))[1].\]
Then $\xi_{H}$ lifts to $F(A\dual)$.
\end{prop}

\begin{pf}
Let $V\subset T[1]$, taking the dual of the exact sequence 
\[ 0\mapor{}V\mapor{}V\somdir{}{}(\symm{}{2}V)\mapor{}
\symm{}{2}V\mapor{}0\]
we get a small extension of finite dimensional nilpotent graded 
algebras
\[ 0\mapor{}(\symm{}{2}V)\dual\mapor{}
V\dual\somdir{}{}(\symm{}{2}V)\dual\mapor{}V\dual 
\mapor{}0.\]
Let 
$ob_{V}\colon F(V\dual)=\Hom_{\cat{G}}(V,T[1])\to
F((\odot^{2}V)\dual[1])=\Hom_{\cat{G}}(\odot^{2}V,T[2])$
be the associated obstruction map and let $q_{V}=ob_{V}(\xi_{V})$. 
Note that $q_{V}$ is, up to sign, the primary obstruction map 
introduced in \ref{primobs}.
By base 
change, for every $V\subset H\subset T[1]$, $q_{V}$ is the obstruction 
to lifting $\xi_{H}$ to $F(H\dual\oplus(\odot^{2}V)\dual)$. 
In particular if the image of $q_{V}(\odot^{2}V)\subset H[1]$, by 
Lemma~\ref{torsorobstruction} the codifferential 
\[\left(\begin{array}{cc}
0&-q_{V}\\ 0&0\end{array}\right)\colon
H\somdir{}{} (\symm{}{2}V)\to (H\somdir{}{}(\symm{}{2}V))[1]\]
is the unique codifferential in $A=H\oplus \odot^{2}V$ such that 
$\xi_{H}$ lifts to $F(A\dual)$. Glueing all the $-q_{V}$, where $V$ 
range over all finite dimensional graded subspaces of $T[1]$, we get 
the required map $Q^{1}_{2}$.\qed\end{pf}

\begin{thm} Let $Q\colon C(T)\to C(T)[1]$ be the linear map
defined by 
\[Q(v_{1}\odot\ldots\odot v_{n})=
\sum_{\sigma\in S(2,n-2)}\epsilon(\sigma)
Q^{1}_{2}(v_{\sigma(1)}\odot v_{\sigma(2)})\odot 
v_{\sigma(3)}\odot\ldots\odot v_{\sigma(n)}.\]
Then $Q$ is a codifferential and then $(C(T),Q)$ is a minimal 
$L_{\oo}$ algebra.\end{thm}

\begin{pf} According to standard formula for $L_{\oo}$ algebras, see 
e.g.
\cite{penkava}, the map $Q$ is a coderivation of degree 1 and $Q^{2}=0$ if 
and only if $Q^{2}(\odot^{3}T[1])=0$.\\
Let $V\subset W\subset Z\subset T[1]$ be finite dimensional subspaces 
such that $Q(\odot^{2}V)\subset W[1]$ and $Q(\odot^{2}W)\subset 
Z[1]$. Then $B=Z\oplus(\odot^{2}W)\oplus(\odot^{3}V)$ 
is a subcoalgebra of $C(T)$ such that $Q(B)\subset B$.
Let $\delta\colon (Z\oplus(\odot^{2}W))\dual\to 
(\odot^{3}V)\dual[2]$ be the transpose of $Q^{2}$, 
by Proposition~\ref{vanishing} we have 
$\delta F((Z\oplus\odot^{2}W)\dual)=0$. Since 
the element $\xi_{Z}$ lifts to some 
$\xi\in F((Z\oplus(\odot^{2}W))\dual)$ and $\delta(\xi)=0$ we also 
have that $\delta=0$ and then also $Q^{2}(B)=0$.\\
The $L_{\oo}$ algebra $(C(T),Q)$ is clearly minimal and, since 
$Q^{1}_{j}=0$ for every $j\ge 3$ it comes from a graded 
Lie algebra $(T,[,])$.\qed\end{pf}

\bigskip
\section{Geometric deformation functors}

In this section we are interested, given a deformation functor $F$, 
to find sufficient conditions for the existence of a 
$L_{\oo}$-algebra $V$ such that $F=\Def_{V}$. 
A straightforward generalization of the transpose of the proof of 
Corollary~\ref{rapprecor} in the set-up of coalgebras shows that, if 
$T^{i}F$ is finite dimensional for every $i$, then there exists a 
structure of minimal $L_{\oo}$-algebra on $T=TF$ such that 
$F=\Def_{T}$.\\ 

Here we consider a special class of functors that, in our 
opinion, contains all the ``interesting'' deformation functors.\\

By extending the analogue definition for classical deformation 
functors \cite{FM1} we shall say that a predeformation functor $G$ is {\em 
homogeneous} if for every pair of morphisms $A\to C$, $B\to C$ in the 
category $\AN$ we have 
\[ G(A\times_{C}B)=G(A)\times_{G(C)}G(B).\]
We say that a deformation functor $F$ if {\em geometric} if $F=G^{+}$ 
for some homogeneous predeformation functor $G$. 
It is not clear to us whether every deformation functor is geometric; 
in any case we strongly suspect that all the concrete example of 
deformation functors having some algebro-geometric interest have this 
property. All the examples of deformation functors considered in 
this paper are geometric.

\begin{thm}\label{struttura}
Let $F$ be a geometric deformation functor and $T=TF$. Then there 
exists a minimal $L_{\oo}$ structure on $T$
and an isomorphism of functors
$\xi\colon \Def_{T}\to F$.
\end{thm}

Before proving the theorem we need some preparatory results.

\begin{defn} Let $(C(V),Q)$ be a $L_{\oo}$ algebra and 
$\{v_{i}\}$ a fixed  homogeneous basis of $V[1]$. We shall say that 
a dg-subcoalgebra $A\subset C(V)$ is {\em standard} if $A$ has a 
basis of monomials $v_{i_{1}}\odot\ldots\odot v_{i_{r}}$; in 
particular every standard dg-subcoalgebra is homogeneous.
\end{defn}

\begin{lem}\label{L-infinitounione} 
Let $(C(V),Q)$ be a minimal $L_{\oo}$-algebra, then 
it is the union of its 
finite dimensional  standard dg-subcoalgebras (with respect to a fixed 
homogeneous basis of $V[1]$).
\end{lem}
\begin{pf}  For every finite dimensional subspace $W\subset C(V)$ 
there exist $n>0$ and a finite dimensional standard subspace $H_{n}$ 
such that $W\subset \oplus_{i=1}^n(\odot^i H_n)$. For 
every sequence of finite dimensional standard subspaces 
\[0\subset H_n\subset H_{n-1}\subset\ldots \subset H_1\subset V[1]\]
we can consider the standard
subcoalgebra $C=\oplus_{i=1}^n(\odot^i H_i)$.
The condition $Q(C)\subset C$ is equivalent to 
\[Q^1_j(\symm{}{j}H_s)\subset H_{s+1-j},\qquad \forall s\ge j\ge 1.\] 
Assume we have the subspaces $H_n,\ldots, H_r$ satisfying the above 
conditions. If $r>1$ then we may take $H_{r-1}$ any finite dimensional 
standard subspace containing 
$\sum_{j\ge 2}Q^1_j(\odot^j H_{r+j-2})$.
\qed\end{pf}

\begin{lem}\label{limiteinverso} 
Let $(C(V),Q)$ be a fixed minimal $L_{\oo}$-algebra and $\{v_{i}\}$ a fixed  
homogeneous basis of $V[1]$;  
denote by $C_n$ the dg-subcoalgebra  $\oplus_{i=1}^n\odot^{i}(V[1])$. 
For a fixed integer $n\ge 2$ let $\sA$ be the set of finite dimensional 
standard dg-subcoalgebras of $C_n$; $\sA$ is a directed set with the 
ordering given by inclusions.  
Let $f\colon \sA\to \cat{Set}$ be a contravariant functor such that:
\begin{enumerate}
\item  if $A\subset C_{n-1}$, then  
$f(A)=\{*\}$ is the one-point set.   
\item $f(A+B)\to f(A)\times_{f(A\cap B)}f(B)$ is bijective.
\item $f(A)\not=\emptyset$ for every  $A\in \sA$.
\end{enumerate}
Then 
\[f(\sA):=\smash{\mathop{liminv}\limits_{A\in \sA}}f(A)\not=\emptyset.\]
\end{lem}

\begin{pf}
We shall say that a subset $\Lambda\subset\sA$ is saturated if:
\begin{enumerate}
\item If $A\subset C_{n-1}$ then $A\in \Lambda$.
\item If $A,B\in\sA$, $A\subset B$ and $B\in\Lambda$ then $A\in\Lambda$.
\item If $A,B\in\Lambda$ then $A+B\in\Lambda$.
\end{enumerate}
For every saturated $\Lambda$ let $f(\Lambda)$ be the set of coherent 
sequences $x\in\prod_{A\in \Lambda}f(A)$.\\
Let $\sB$ the set of all pairs $(\Lambda,x)$ with $\Lambda\subset\sA$ 
saturated, $f(\Lambda)\not=\emptyset$ 
and $x\in f(\Lambda)$. On $\sB$ there exists a natural ordering: 
$(\Lambda_1,x_1)\le (\Lambda_2,x_2)$ if and only if $\Lambda_1\subset 
\Lambda_2$ and $x_2$ extends $x_1$. 
Clearly $\sB\not=\emptyset$ and 
every chain in $\sB$ has an upper bound. By Zorn's lemma 
there exists a maximal element $(\Lambda,x)\in \sB$; we want to prove 
that $\Lambda=\sA$.\\   
Assume $\Lambda\not=\sA$ and let $A\in \sA-\Lambda$ such that the 
dimension of $A$ is minimum. Necessarily we must have 
$\dim_{\K}(A)=\dim_{\K}(A\cap C_{n-1})+1$, $A\cap B\subset 
C_{n-1}$ and $f(A+B)=f(A)\times f(B)$ for every $B\in\Lambda$.
Let $\Lambda'$ be the union of $\Lambda$ and $\{A+B\,|\, 
B\in\Lambda\}$.
Choose an element $x_{A}\in f(A)$, then it is uniquely determined an 
element $x_{A+B}\in f(A+B)$ for every $B\in \Lambda$. 
We need to prove that, if $C\subset A+B$ then either $C=A+B$ or $C\in\Lambda'$ 
  and $x_{A+B}$ extends $x_{C}$.\\
Assume $C\not=A+B$, since $A,B,C$ are standard we have $C=(A\cap 
C)+(B\cap C)$. If $A\subset C$ then $C=A+(B\cap C)\in\Lambda'$, since 
$x_{B}$ extends $x_{B\cap C}$ also $x_{A+B}$ extends $x_{C}$. If 
$A\cap C\not=A$ then $C\cap A\in\Lambda$, $C\cap A\subset C_{n-1}$ and
$f(C)=f(C\cap B)$; since $x_{A+B}$ extends $x_{B}$ it also extends $x_{C}$.
\qed\end{pf}

\begin{rem} It is clear that Lemma~\ref{limiteinverso} 
also works if the 
codifferential $Q$ is defined only in the subcoalgebra $C_{n}$.
\end{rem}

\begin{pfof}{\ref{struttura}}
Denote by $S_n=\oplus_{i=1}^n(\odot^i(T[1]))$. We recall 
that for every graded vector space 
$I\in \AN\cap\cat{G}$  we have
\[F(I)=\Hom_{\cat{G}}(I\dual,T[1])=\Hom_{\cat{DGC}}(I\dual,S_1),\]
where $S_1$ is endowed with the codifferential $d_1=0$. In particular for 
every finite dimensional standard subspace $J\subset T[1]$ we have 
$F(J\dual)=\Hom_{\cat{G}}(J,T[1])$: we denote by $\xi_J\in F(J\dual)$ the 
element corresponding to the inclusion $J\hookrightarrow T[1]$ and by 
\[ \xi_1=(\xi_J)\in\prod_{J\subset T[1]}F(J\dual).\]
Let $G$ be a homogeneous predeformation functor such that $F=G^{+}$ 
and let $\pi\colon G\to F$ the projection.\\ 
We now construct recursively a sequence $(d_n,\eta_n)$ such that for every 
$n\ge 1$
\begin{description}
\item{$A_n$:~} $d_n\colon S_n\to S_{n-1}[1]\subset S_{n}[1]$ 
is a codifferential extending $d_{n-1}$.
\item{$B_n$:~} $\eta_n\in \mathop{liminv}F(G\dual)$, where: 
\begin{enumerate} 
\item the limit is 
taken over all the finite dimensional standard dg-subcoalgebras of 
$(S_n,d_n)$.
\item $\eta_{1}$ is a lifting of $\xi_{1}$.
\item $\eta_n$ is a lifting of $\eta_{n-1}$ if $n\ge 2$.
\end{enumerate}\end{description}     

For every standard subspace $J\subset T[1]$ let $f(J)\subset 
G(J\dual)$ be the set of liftings of $\xi_{J}$; the functor $f$ 
satisfies the hypothesis of Lemma~\ref{limiteinverso} and then there 
exists $\eta_{1}$ that lifts $\xi_{1}$.
Let $n\ge 2$ be a fixed integer and assume are given $(d_1,\eta_1),\ldots,
(d_{n-1},\eta_{n-1})$ satisfying  $A_i$, $B_i$, $i=1,\ldots,n-1$.\\
Let $A\subset S_n$ be a fixed standard finite dimensional subcoalgebra and 
denote $B=A\cap S_{n-1}$, $J=A\cap S_1$. 
Possibly enlarging $B$, we may assume  that 
$d_{n-1}(B)\subset B$ and $d_{n-1}$ can be extended to a coderivation 
on $A$; dualizing we get a small extension of finite 
dimensional graded nilpotent algebras
\[0\mapor{}B^{\perp}\mapor{}A\dual\mapor{}B\dual\mapor{}0,\]
such that $B^{\perp}$ and $B\dual$ have respective 
differential equal to $0$ and ${}^td_{n-1}$.\\ 
Let $\xi_{B}=\pi(\eta_{B})\in F(B\dual)$ be the projection of 
$\xi_{n-1}=\pi(\eta_{n-1})$. Since $\xi_{B}$ lifts $\xi_{J}$
and 
$J\dual=B\dual/(B\dual)^2$ 
the map 
\[ [B\dual,B^{\perp}[2]]\to F(B^{\perp}[2]),\qquad
[f]\to f(\xi_{B})\]
is injective and, by Proposition~\ref{vanishing}, there exists a differential 
on $A\dual$ making the exact sequence 
$0\to B^{\perp}\to A\dual\to B\to 0$ a small extension in 
$\AN$.\\
Let $o\in F(B^{\perp}[1])=\Hom_{\cat{G}}(A/B,T[2])$ be the 
obstruction to the lifting of $\xi_{B}$ to $F(A\dual)$; 
possibly enlarging $B$, we may suppose that 
$o\in\Hom_{\cat{G}}(A/B,J[1])$ and then there exists an 
unique morphism $\phi\colon B\dual\to B^{\perp}[1]$ such that 
$o=\phi(\xi_{B})$. By Lemma~\ref{torsorobstruction} there exists an 
unique differential on $A$ making  
$0\to B^{\perp}\to A\dual\to B\to 0$ a small extension in 
$\AN$ and such that $\xi_{B}$ lifts to $F(A\dual)$. The transpose of 
this differential is a codifferential in $A$.\\
Clearly all these codifferentials in the finite dimensional 
homogeneous subcoalgebras of 
$S_n$ can be glued together and we get a global codifferential 
$d_n\colon S_n\to S_{n-1}[1]$.\\
For every $A\subset S_n$ finite dimensional standard
dg-subcoalgebra let
\[ f(A)=\{x\in G(A\dual)\,|\, 
x\hbox{ is a lifting of } \eta_{A\cap S_{n-1}}\}.\]
By Lemma~\ref{smoothpredef} and
the choice of $d_n$, $f(A)\not=\emptyset$ for every $A$. The functor 
$f$ satisfies the hypothesis of Lemma~\ref{limiteinverso}
and then $\mathop{liminv}f(A)\not=\emptyset$. 
By construction every 
$\eta_n\in \mathop{liminv}f(A)$ satisfies the condition $B_n$.\\
The sequence $\eta=(\eta_n)$ gives a natural transformation 
$\xi\colon MC_T\to G$ by the rule:
\[ \Hom_{\cat{DGC}}(A\dual,C(T))\ni \alpha\to 
\eta(f)={}^t\alpha(\eta_D)\in F(A),\]
where $D$ is a finite dimensional standard dg-subcoalgebra of $C(T)$ 
such that $\alpha(A\dual)\subset D$.\\
By construction the induced morphism $\Def_T=MC_T^+\to F$ is an 
isomorphism on tangent spaces and then, according to 
Corollary~\ref{CorIMT} it is an isomorphism of functors.
\qed\end{pfof}

%\section{biblio}

\begin{tabular}{l} 
Marco Manetti\\
 Dipartimento di Matematica ``G. Castelnuovo'',\\
 Universit\`a di Roma ``La Sapienza'',\\ 
 Piazzale Aldo Moro 5, I-00185 Roma, Italy.\\
 manetti@mat.uniroma1.it,~~~~~
http://www.mat.uniroma1.it/people/manetti/\\
 \end{tabular}

					\end{document}